\documentclass[10pt,reqno]{amsart}
\usepackage{hyperref}
\usepackage{tikz}
\usetikzlibrary{decorations.pathreplacing}

\usepackage{tikz}
\usetikzlibrary{patterns.meta,
                decorations.pathreplacing,
                    calligraphy,
                positioning
                }
\usepackage{soul}
\usepackage{amssymb}
\numberwithin{equation}{section}
\usepackage{hyperref}
\usepackage{float}
\usepackage{pgf,tikz}

\newtheorem{theorem}{Theorem}[section]
\newtheorem{definition}[theorem]{Definition}

\newtheorem{lemma}[theorem]{Lemma}

\newtheorem{example}[theorem]{Example}

\newtheorem{remark}[theorem]{Remark}

\newtheorem{corollary}[theorem]{Corollary}

\newtheorem*{notation*}{Notation}

\begin{document}

\title{The first syzygy of Hibi rings associated with planar distributive lattices}

\author[P. Das]{Priya Das}
\address{Indian Institute of science education and research kolkata, West
Bengal, India}
\email{priya.math88@gmail.com}
\author[H.~Mukherjee]{Himadri Mukherjee}
\address{BITS Pilani K.K. Birla Goa campus, India}
\email{himadrim@goa.bits-pilani.ac.in}

\thanks{AMS Classification 2010. Primary: 05E40, 06D75, 03G10, 16E05}
\maketitle
\begin{abstract}
  In this article, we give explicit minimal generators of
the first syzygy of the Hibi ring for a planar distributive lattice in terms of sublattices. We also give a characterization when it is linearly related and derive an exact formula for the first Betti number of a planar distributive lattice.
\end{abstract}

\section{Introduction}
Let $\mathcal{L}$ be a finite distributive lattice and let $S=K[x_{\alpha}: \alpha \in \mathcal{L}]$ be a polynomial ring over a field 
$K$. For $\alpha,\beta \in \mathcal{L}$,
$f_{(\alpha,\beta)}=x_{\alpha}x_{\beta}- x_{\alpha \vee \beta}x_{\alpha \wedge \beta} \in S$, is called a diamond relation
if $\alpha \nsim \beta$ (we say, $\alpha$ is not comparable to $\beta$). Let $I=\langle f_{(\alpha,\beta)}: \alpha \nsim \beta, \alpha,\beta \in \mathcal{L}\rangle \subset S$ be the associated ideal, called a Hibi ideal. The associated ring $R[\mathcal{L}]=S/I$ is called a Hibi ring. 

Hibi rings were introduced by Hibi in his paper \cite{hibi00}, he studied $R[\mathcal{L}]$
as an algebra with straightening laws on $\mathcal{L}$ over the field $K$. He demonstrated that $R[\mathcal{L}]$ is Gorenstein when the set of join irreducible elements
$J=\{ z \in \mathcal{L}: x \vee y=z \Rightarrow x=z$ or $y=z \}$ is pure i.e. all maximal chains have the same length. In the same paper, it is shown $R[\mathcal{L}]$ is a Cohen-Macaulay normal domain if $\mathcal{L}$ is distributive. Further, Hibi ideal is prime if and only if $\mathcal{L}$ is distributive \cite{hibi00}.

  Herzog-Hibi (\cite{herzog_monomial}, Theorem 10.1.3), showed that
the diamond
relations $f_{(\alpha,\beta)}$ form a Gr\"{o}bner basis with respect to the reverse lexicographic term order on $S$ extending the order of $\mathcal{L}$.
The ring appears in  geometric context in Lakshmibai-Gonciulea \cite{gonciulea_singular}, where it is proved that
the Schubert varieties in $G_{d,n}$ degenerate flatly to the toric varieties
 $X(I_{d,n})(=X_{d,n}$ as in \cite{gonciulea_singular}) which are varieties associated to the lattices $I_{d,n}=\{
 (i_1,\ldots,i_d): 1 \leq i_1 < \ldots < i_d \leq n\}$. In the same paper, using the degeneration, several geometric properties for the
 Schubert varieties were derived. The
 question of singularity of the variety associated with the Hibi ideal is dealt with in Wagner, \cite{wagner_singularities}.  Lakshmibai-Mukherjee, \cite{himadri}, gives a standard monomial basis for the tangent space at a singular point for these varieties, they used the standard monomial basis to compute the tangent cone at these singular point, \cite{himadri_thesis}. The question of free resolution of Hibi rings has been addressed in a few recent papers, \cite{ene_syzygies}, \cite{hibi_reg}, for a  modular non-distributive lattice. Rodica-Ene-Hibi give a lower bound for the regularity of $I$ (Theorem 2.2, \cite{hibi_reg}) and in (Theorem 1.1, \cite{herzog} ) authors have given the regularity of Hibi ring of any finite distributive lattice in terms of its poset of join-irreducibles.
 Hibi rings over a finite distributive lattice are further studied by several authors in, \cite{aramova}, \cite{brown_singular}, \cite{brown-lucky}, \cite{ene_join} to mention but a few. 
 Construction of an explicit minimal graded free resolutions of classes of graded modules is difficult \cite{peeva_graded} i.e. to find the explicit generators of the modules is a challenging problem. In particular for distributive lattices, the following problems may be asked 
\begin{enumerate}
  \item What are the explicit minimal generators of the first syzygy of $R[\mathcal{L}]$, in terms of the the lattice $\mathcal{L}$? 
  \item Find an exact formula for the first Betti number of $R[\mathcal{L}]$.
 \end{enumerate}
 The present paper deals with these two questions for finite planar distributive lattices. Our approach is that the generators are linked with certain kind of sublattices of the lattice $\mathcal{L}$, as a result it is easier to verify conditions such as linearity of the generators.

In our theorem (Theorem \ref{main} and Theorem \ref{last}) we have totally characterized the lattices for which the syzygy is linear, although such theorem is available in literature \cite{ene_syzygies} our result has a different description that relies on the number of join meet irreducibles. In the very useful survey paper (which also contains a few new results by the author) \cite{ene_syzygies}, Ene answers when $R[\mathcal{L}]$ has a linear resolution and characterizes the lattices for which $I$ is linearly related. 

The following is our main theorem that totally characterizes the syzygy generators for the Hibi ideal of a planar distributive lattice in terms of certain relations arizing out of sublattices. The novelty of this theorem is that a-priori there is no reason why the minimal generators of the syzygy module should even be related to certain sublattices, let alone a small number of such sublattices. 

Sublattices isomorphic to the following three, we call them classifying sublattices, in the figure \ref{figure} give rise to all the linear first syzygy elements. The first one (a) gives rise two two types of generators we call them $S_1,S_2$ the next one (b) gives rise to the type $L$ and the last one (c) gives rise to two generators $B_1,B_2$. 

\begin{figure}[H] 
 \begin{tikzpicture}[scale=.2]
 \draw (0,0)--(-2,2)--(0,4)--(2,6)--(4,4)--(2,2)--(0,0);
 \draw (2,2)--(0,4);
 \draw (12,0)--(10,2)--(8,4)--(10,6)--(12,8)--(14,6)--(12,4)--(14,2)--(12,0);
 \draw (10,2)--(12,4)--(10,6);
 \draw (22,0)--(20,2)--(18,4)--(20,6)--(22,8)--(24,6)--(26,4)--(24,2)--(22,0);
 \draw (20,2)--(22,4)--(24,6);
 \draw (24,2)--(22,4)--(20,6);
 \draw (-.5,-1) node[anchor=north west] {$(a)$};
 \draw (11.5,-1) node[anchor=north west] {$(b)$};
 \draw (21.5,-1) node[anchor=north west] {$(c)$};
 \begin{scriptsize}
 \fill [color=black] (0,0) circle (1.5pt);
\fill [color=black] (-2,2) circle (1.5pt);
\fill [color=black] (0,4) circle (1.5pt);
\fill [color=black] (2,6) circle (1.5pt);
\fill [color=black] (4,4) circle (1.5pt);
\fill [color=black] (2,2) circle (1.5pt);
\fill [color=black] (12,0) circle (1.5pt);
\fill [color=black] (10,2) circle (1.5pt);
\fill [color=black] (8,4) circle (1.5pt);
\fill [color=black] (10,6) circle (1.5pt);
\fill [color=black] (12,8) circle (1.5pt);
\fill [color=black] (14,6) circle (1.5pt);
\fill [color=black] (12,4) circle (1.5pt);
\fill [color=black] (14,2) circle (1.5pt);
\fill [color=black] (22,0) circle (1.5pt);
\fill [color=black] (20,2) circle (1.5pt);
\fill [color=black] (18,4) circle (1.5pt);
\fill [color=black] (20,6) circle (1.5pt);
\fill [color=black] (22,8) circle (1.5pt);
\fill [color=black] (22,4) circle (1.5pt);
\fill [color=black] (24,6) circle (1.5pt);
\fill [color=black] (26,4) circle (1.5pt);
\fill [color=black] (24,2) circle (1.5pt);
 \end{scriptsize}
\end{tikzpicture}
\label{figure}
\caption{Classifying sublattices}
\end{figure}
\begin{theorem}[Theorem \ref{syz}]\label{main_1}
 Syz$^{1}_{1} R[\mathcal{L}]$ is generated by the types $S_1, S_2, L, B_1, B_2$.
\end{theorem}
In Theorem \ref{main_1}, Syz$^{1}_{1}R[\mathcal{L}]$ denotes linear generators of first syzygy of $R[\mathcal{L}]$.
This main result of this paper is quoted in the theorem below. This theorem along with the quadratic elements coming from the \textit{non-expressible} generators do form a linearly independent set of generators for the first syzygy module for the Hibi ideals, quoted in the Theorem \ref{main_2}. In particular, it also characterizes when $R[\mathcal{L}]$ is linearly related.

\begin{theorem}[Theorem \ref{main}]\label{main_2}
 Syz$^{2}_{1} R[\mathcal{L}]$ is trivial if and only if every pair of comparable diamonds in $\mathcal{D}_2$ are expressible.
\end{theorem}
In Theorem \ref{main_2}, Syz$^{2}_{1} R[\mathcal{L}]$ denote the set of all diamond type (for definition see \ref{Diamond}) generators of first syzygy of $R[\mathcal{L}]$.

 Let $JM$ be the set of join-meet irreducible elements of a lattice $\mathcal{L}$ and let us denote $k(\mathcal{L})=\# \{(\theta_i,\theta_j): \theta_i \nsim \theta_j, \theta_i,\theta_j \in JM\}$.

\begin{theorem}[Theorem \ref{last}]
  With the above notations, we have the following
\begin{enumerate}
 \item When $k(\mathcal{L})=1$, the first syzygy of $R[\mathcal{L}]$ is linear.
 \item When $k(\mathcal{L})=2$,
 \begin{enumerate}
  \item If the lattice is Figure 7, then first syzygy is non linear.
  \item  Else if the lattice is Figure 8, then first syzygy is linear if and only if $ht(\theta_2 \wedge \theta_3)
 -ht(\theta_1 \wedge \theta_2)=1$ or $ht(\theta_2 \vee \theta_3)-ht(\theta_1 \vee \theta_2)=1$, where $\theta_1,\theta_2,
 \theta_3 \in JM$.
 \end{enumerate}
 \item If $k(\mathcal{L}) \geq 3$ then the first syzygy is non linear.
\end{enumerate} 
\end{theorem}

The paper is structured as follows. In Section \ref{pre}, we collect basic notations, terminology, and results that will be used in the paper. The first syzygy of Hibi rings is discussed in Section \ref{Hibi_ring}. Explicit expression for the first Betti number for planar distributive lattices has been discussed in Section \ref{Betti}.

\section{Preliminaries}\label{pre}
In this section, we introduce definitions and results needed in the next section. Let $(P,\leq)$ be a poset with partial order relation $\leq$. The \textit{height} of an element $a\in P$, denoted by $ht(a)$ is the supremum of length of chains $C$, descending from $a$. A \textit{lattice} $\mathcal{L}$ is a poset such that for all $x,y \in \mathcal{L}, x\vee y, x\wedge y$ exists. A lattice $\mathcal{L}_1$ is called \textit{sublattice} of the lattice $\mathcal{L}$ if $\mathcal{L}_1\subset \mathcal{L}$ and $x,y \in \mathcal{L}_1$ implies $x\vee y \in \mathcal{L}_1$ and $x\wedge y \in \mathcal{L}_1$. Let $C_{m+1}$ be the chain $1<2<\ldots <m+1$, then  $G(m,n)=C_{m+1}\times C_{n+1}$ be the product lattice, we call $G(m,n)$ an $m\times n$ \textit{grid lattice}. A finite lattice $\mathcal{L}$ is called \textit{planar} if and only if there is an embedding of $\mathcal{L}$ in grid lattice $G(m,n)$ for some $m,n$. In a lattice $\mathcal{L}$, an element $x \in \mathcal{L}$ is said to be \textit{join irreducible} if $x=y\vee z$ then $x=y$ or $x=z$ for $x,y,z \in \mathcal{L}$ and is called \textit{meet irreducible} if $x=y\wedge z$ then $x=y$ or $x=z$.

A non-empty subset $\mathcal{I}(\mathcal{L})$ of a poset $P$ is called \textit{poset ideal} if for any $x\in\mathcal{I}(\mathcal{L})$ and $y\in P$ if $y\leq x$ then $y \in \mathcal{I}(\mathcal{L})$. Now we state the fundamental theorem of distributive lattice which is implicitly used in the paper. \begin{theorem}[Birkhoff Theorem, \cite{gratzer_lattice}]Let $\mathcal{L}$ be a distributive lattice and $P$ be the poset of join-irreducible elements. Then $\mathcal{L}\simeq I(P)$.
\end{theorem}
Let $S=K[x_1,x_2,\ldots,x_n]$ be a polynomial ring over a field $K$. If $>$ is a monomial order then for any $f\in S$, the initial term of $f$, \textit{in$_>(f)$}, is the greatest term of $f$ with respect to the order $>$. We will use, in$(f)$, when the order is understood. Let $m,n$ be two monomials of the polynomial ring $S$, let $a=(a_1,\ldots,a_n), b=(b_1,\ldots,b_n)$ the multi indices of $m,n$ respectively. Then the \textit{reverse lexicographic order} is defined by $m>_{rlex} n$ if and only if deg$(m)$ $>$ deg$(n)$ or deg $(m)$=deg $(n)$ and $a_i<b_i$ for the least integer $i$ with $a_i\neq b_i$.

To find the first syzygy we are frequently using the following (\cite{eisenbud_commutative}, for more details).
Let $F$ be a free module with basis and $M$ be a submodule of $F$ generated by monomials $m_1,\ldots,m_t$. Let $\phi: \oplus_{j=1}^{t} S \varepsilon_{j} \longrightarrow F$ 
be a homomorphism defined by $\varepsilon_i \longmapsto m_j$, whose image is $M$. For each pair of indices $i$ and $j$ such that $m_i, m_j$ involve same basis element of $F$, we define $m_{ij}=m_i/gcd(m_i,m_j)$ and let ker$\phi=\sigma_{ij}$ is given by $\sigma_{ij}=m_{ji}\varepsilon_i-m_{ij}\varepsilon_j$.
\begin{lemma}[Lemma \cite{eisenbud_commutative}] With the same notation as above, ker$\phi$ is generated by the $\sigma_{ij}$. 
\end{lemma}

\section{First syzygy of Hibi rings}\label{Hibi_ring}
In this section, we find explicitly the generators of the first syzygy of Hibi ring for a planar distributive lattice and also find a condition for which the first syzygy linear.

Let $f_1, f_2, \ldots, f_n$ be all the diamond relations in $I$. Let us denote $F=\oplus_{i=1}^{n} S(-2)$ 
and $g_1,g_2,\ldots,g_n$ be the generators of $F$,
 $\phi: F \longrightarrow S$, $g_i \longmapsto f_i$, ker$(\phi)$ is the
first syzygy of $R[\mathcal{L}]$. Let us denote ker$(\phi)$=Syz$^{1}R[\mathcal{L}]$ and Syz$^{1}_{1} R[\mathcal{L}]$ denote linear 
generators of the  $R[\mathcal{L}]$-module Syz$^{1}R[\mathcal{L}]$.\\ \\
Let $Q=\{\alpha,\beta,\gamma,\delta\}$ be a poset $\leq$ defined by $\alpha \leq \beta,\gamma,\delta$, where $\beta,\gamma,\delta$ are mutually non-comparable. Let the lattice of the poset 
ideals of $Q$ be called a cube lattice.

\begin{lemma}\label{planar} $\mathcal{L}$ is a planar distributive lattice if and only if every join irreducible (meet irreducible) 
 $\beta \in \mathcal{L}$ is covered (covers) by at most two join irreducibles (meet irreducibles) in $\mathcal{L}$.
\end{lemma}
\begin{proof}
 $\Rightarrow:$ Let $\mathcal{L}$ be a planar distributive lattice and $\beta \in \mathcal{L}$ be a join-irreducible element. Let us assume that it is covered by three join-irreducible elements $x,y,z$ in $\mathcal{L}$. Then 
$x \vee y, x \vee z, y\vee z$ are non-comparable. Hence there will be a cube sublattice in $\mathcal{L}$, which is 
non-planar.
\\
$\Leftarrow$: Let every join-irreducible $\beta$ in $\mathcal{L}$ covered by at most two join-irreducible in $\mathcal{L}$, let this
property be called (*). Then 
every sublattice of $\mathcal{L}$ also has the property (*). Now if $\mathcal{L}$ is not planar then there exist a cube 
sublattice in $\mathcal{L}$, then this sublattice violates the property (*), as the cube has three join-irreducible elements
covering the minimal element which is a join-irreducible. Hence we arrive at a contradiction. By duality we get the result with 
meet irreducibles.
\end{proof}
 \begin{lemma}\label{com}$\mathcal{L}$ is a planar distributive lattice if and only if for every $x,y,z \in \mathcal{L}, x\vee y, y \vee z, x \vee z$ are comparable.  
 \end{lemma}
 \begin{proof}  $\Rightarrow$: Let $\mathcal{L}$ be a planar distributive lattice. Let $x,y,z \in \mathcal{L}$, then 
the sublattice generated by $x \wedge y, x \wedge z, y \wedge z$ are also planar. If $x \wedge y, x \wedge z, y \wedge z$
are not comparable to each other then we will get three join-irreducibles covering minimal element. Therefore by Lemma \ref{planar},
we will arrive at a contradiction. Hence $x\wedge y \sim x \wedge z$ and by duality other also holds.
\\
$\Leftarrow$: Let the given conditions hold.
 Then we have the cases $(i)\, x\vee y \leq x \vee z,\, x\wedge y \leq x \wedge z\, (ii)\, x \vee y \geq x \vee y, \,
 x \wedge y \leq x \wedge z\,
 (iii)\,  x \vee y \geq x \vee y,\, x \wedge y \leq x \wedge z$. Then the only possibilities to get the sublattices as in Figure 1 and Figure 2, which are all planar sublattices. Hence $\mathcal{L}$ is planar.
 \end{proof} 
 Let us take two diamond relations $f_{(\alpha_1,\beta_1)}=x_{\alpha_1}x_{\beta_1}-x_{\alpha_1 \vee \beta_1}
x_{\alpha_1 \wedge \beta_1}$ and $f_{(\alpha_2,\beta_2)}=x_{\alpha_2}x_{\beta_2}-x_{\alpha_2 \vee \beta_2}
x_{\alpha_2 \wedge \beta_2}$ coming from the diamonds $(\alpha_1,\beta_1)$ and $(\alpha_2,\beta_2)$ respectively. Then
with respect to reverse lexicographic order $>$, we have in$(f_{(\alpha_1,\beta_1)})=x_{\alpha_1}x_{\beta_1}$,
in$(f_{(\alpha_2,\beta_2)})=x_{\alpha_2}x_{\beta_2}$. Then we have the following cases
\vskip 1mm \noindent
\textsc{Case 1:} $\alpha_1=\alpha_2$ or $\alpha_1=\beta_2$ or $\beta_1=\alpha_2$ or $\beta_1=\beta_2$.
\vskip 1mm \noindent
\textsc{Case 2:} $\alpha_1,\alpha_2,\beta_1,\beta_2$ all are different (it is dealt with definition \ref{Diamond} and after).

 For \textsc{Case 1}, without a loss of generality, let us suppose $\alpha_1=\alpha_2$. Then we have 
in$(f_{(\alpha_1,\beta_1)})=x_{\alpha_1}x_{\beta_1}$, in$(f_{(\alpha_1,\beta_2)})=x_{\alpha_1}x_{\beta_2}$ and therefore,
Lcm$(f_{(\alpha_1,\beta_1)},f_{(\alpha_1,\beta_2)})=x_{\alpha_1}x_{\beta_1}x_{\beta_2}$. So,
\begin{align*}
 S(f_{(\alpha_1,\beta_1)},f_{(\alpha_1,\beta_2)})=&x_{\beta_2}f_{(\alpha_1,\beta_1)}- x_{\beta_1}f_{(\alpha_1,\beta_2)}\\
 =& x_{\beta_1}x_{\alpha_1 \vee \beta_2}x_{\alpha_1 \wedge \beta_2}- x_{\beta_2}x_{\alpha_1 \vee \beta_1}
 x_{\alpha_1 \wedge \beta_1}
\end{align*}
Then either $\beta_1 \sim \beta_2$, let us say $\beta_1 \leq \beta_2$ or $\beta_1 \nsim \beta_2$. So we have the following subcases
\vskip 1mm \noindent
\textsc{Subcase 1:} $\beta_1 \leq \beta_2$
\vskip 1mm \noindent
\textsc{Subcase 2:} $\beta_1 \nsim \beta_2$

\noindent Also by Lemma \ref{com}, for 
$\alpha_1,\beta_1,\beta_2$ we have $\alpha_1 \wedge \beta_1 \sim \alpha_1 \wedge \beta_2$ and 
$\alpha_1 \vee \beta_1 \sim \alpha_1 \vee \beta_2$. Then we have the only following  possibilities
\begin{enumerate}
 \item $\alpha_1 \wedge \beta_1 \leq \alpha_1 \wedge \beta_2$, $\alpha_1 \vee \beta_1 \leq \alpha_1 \vee \beta_2$. 
 \item $\alpha_1 \wedge \beta_1 \geq \alpha_1 \wedge \beta_2$, $\alpha_1 \vee \beta_1 \leq \alpha_1 \vee \beta_2$.
 \item $\alpha_1 \wedge \beta_1 \leq \alpha_1 \wedge \beta_2$, $\alpha_1 \vee \beta_1 \geq \alpha_1 \vee \beta_2$.
\end{enumerate}
Now we verify all the above cases
\vskip 1mm \noindent
\textsc{Subcase 1:} When $\beta_1 \leq \beta_2$.
 \begin{itemize}
 \item From the condition $(1)$, we get a sublattice isomorphic to that in Figure 2. $(a)$.
 
\begin{figure}[H] 
 \begin{tikzpicture}[scale=.3]
 \draw (0,0)--(-2,2)--(-4,4)--(-2,6)--(0,8)--(2,6)--(0,4)--(2,2)--(0,0);
 \draw (-2,2)--(0,4);
 \draw (0,4)--(-2,6);
 \draw (10,2)--(8,4)--(10,6)--(12,8)--(14,6)--(12,4)--(10,2);
 \draw (10,6)--(12,4);
 \draw (22,2)--(20,4)--(18,6)--(20,8)--(22,6)--(24,4)--(22,2);
 \draw (20,4)--(22,6);
 \draw (-1,-1.5) node[anchor=north west] {$(a)$};
 \draw (9,.5) node[anchor=north west] {$(b)$};
 \draw (21,.5) node[anchor=north west] {$(c)$};
 \begin{scriptsize}
 \fill [color=black] (0,0) circle (1.5pt);
\draw[color=black] (0,-.5) node {$\alpha_1\wedge \beta_1$};
\fill [color=black] (-2,2) circle (1.5pt);
\draw[color=black] (-3.5,2) node {$\alpha_1\wedge \beta_2$};
\fill [color=black] (-4,4) circle (1.5pt);
\draw[color=black] (-4.5,4) node {$\alpha_1$};
\fill [color=black] (0,4) circle (1.5pt);
\draw[color=black] (3,4) node {$\beta_1 \vee (\alpha_1 \wedge \beta_2)$};
\fill [color=black] (2,2) circle (1.5pt);
\draw[color=black] (2.7,2) node {$\beta_1$};
\fill [color=black] (-2,6) circle (1.5pt);
\draw[color=black] (-3.7,6) node {$\alpha_1 \vee \beta_1$};
\fill [color=black] (0,8) circle (1.5pt);
\draw[color=black] (0,8.5) node {$\alpha_1 \vee \beta_2$};
\fill [color=black] (2,6) circle (1.5pt);
\draw[color=black] (2.7,6) node {$\beta_2$};
\fill [color=black] (10,2) circle (1.5pt);
\draw[color=black] (10.5,1.5) node {$\alpha_1\wedge \beta_1(=\alpha_1 \wedge \beta_2)$};
\fill [color=black] (8,4) circle (1.5pt);
\draw[color=black] (7.5,4) node {$\alpha_1$};
\fill [color=black] (10,6) circle (1.5pt);
\draw[color=black] (8.5,6) node {$\alpha_1 \vee \beta_1$};
\fill [color=black] (12,8) circle (1.5pt);
\draw[color=black] (12,8.5) node {$\alpha_1 \vee \beta_2$};
\fill [color=black] (14,6) circle (1.5pt);
\draw[color=black] (14.5,6) node {$\beta_2$};
\fill [color=black] (12,4) circle (1.5pt);
\draw[color=black] (12.9,4) node {$\beta_1$};
\fill [color=black] (22,2) circle (1.5pt);
\draw[color=black] (22,1.5) node {$\alpha_1 \wedge \beta_1$};
\fill [color=black] (20,4) circle (1.5pt);
\draw[color=black] (18.5,4) node {$\alpha_1 \wedge \beta_2$};
\fill [color=black] (18,6) circle (1.5pt);
\draw[color=black] (17.5,6) node {$\alpha_1$};
\fill [color=black] (20,8) circle (1.5pt);
\draw[color=black] (21,8.5) node {$\alpha_1 \vee \beta_1(=\alpha_1 \vee \beta_2)$};
\fill [color=black] (22,6) circle (1.5pt);
\draw[color=black] (22.7,6) node {$\beta_2$};
\fill [color=black] (24,4) circle (1.5pt);
\draw[color=black] (24.7,4) node {$\beta_1$};
 \end{scriptsize}
\end{tikzpicture}
\label{figure}
\caption{}
\end{figure}

 \item We get from the condition $(2)$, a sublattice isomorphic to that in Figure 2. $(b)$.
  
  \item  The condition $(3)$, leads us a sublattice isomorphic to that in Figure 2. $(c)$.

  \end{itemize}
\vskip 1mm \noindent
\textsc{Subcase 2:} When $\beta_1 \nsim \beta_2$.
  \begin{itemize}
   \item From the condition $(1)$, we will not get any lattice which satisfies the given condition, since from this we get $\beta_1 \leq \beta_2$, then it will give us a contradiction.
   \item From the condition $(2)$ we get the sublattice isomorphic to that in Figure 3.
   \begin{figure}[H] 
 \begin{tikzpicture}[scale=.3]
 \draw (0,0)--(-2,2)--(-4,4)--(-2,6)--(0,8)--(2,6)--(0,4)--(2,2)--(0,0);
 \draw (-2,2)--(0,4);
 \draw (0,4)--(-2,6);
 \draw (2,6)--(4,4)--(2,2);

 \begin{scriptsize}
  \fill [color=black] (0,0) circle (1.5pt);
\draw[color=black] (0,-.5) node {$\alpha_1\wedge \beta_2$};
\fill [color=black] (-2,2) circle (1.5pt);
\draw[color=black] (-3.5,2) node {$\alpha_1\wedge \beta_1$};
\fill [color=black] (-4,4) circle (1.5pt);
\draw[color=black] (-4.5,4) node {$\alpha_1$};
\fill [color=black] (0,4) circle (1.5pt);
\draw[color=black] (1,4) node {$\beta_1$};
\fill [color=black] (2,2) circle (1.5pt);
\draw[color=black] (3.7,2) node {$\beta_1 \wedge \beta_2$};
\fill [color=black] (-2,6) circle (1.5pt);
\draw[color=black] (-3.7,6) node {$\alpha_1 \vee \beta_1$};
\fill [color=black] (0,8) circle (1.5pt);
\draw[color=black] (0,8.3) node {$\alpha_1 \vee \beta_2$};
\fill [color=black] (2,6) circle (1.5pt);
\draw[color=black] (3.7,6) node {$\beta_1 \vee \beta_2$};
\fill [color=black] (4,4) circle (1.5pt);
\draw[color=black] (4.7,4) node {$\beta_2$};
 \end{scriptsize}
 \end{tikzpicture}
 \label{box}
 \caption{}
 
   \end{figure}
  
 \item The condition $(3)$, will give the same sublattice isomorphic to that in Figure 3, by replacing $\beta_1$ and $\beta_2$ by $\beta_2$ and $\beta_1$ respectively.    
  \end{itemize}

Let us denote 
 $\mathcal{D}_2=\{(\alpha,\beta): \alpha \nsim \beta, \alpha, \beta \in \mathcal{L}\}$. 
   Now using the conditions of the above cases (for  which lattice exists), we prove that the following are the types of the 
 generators of the first syzygy of the Hibi ring.
 
  \begin{lemma}\label{S} $(\textbf{Strip-type})$ Let $(\alpha_1,\beta_1),(\alpha_1,\beta_2) \in \mathcal{D}_2$, for every sublattice 
 isomorphic to Figure 2.(b), with conditions $\alpha_1 \vee \beta_1 \neq \alpha_1 \vee \beta_2, \alpha_1 \wedge \beta_1=
 \alpha_1 \wedge \beta_2, \beta_1 \leq \beta_2, \beta_1 \leq \alpha_1 \vee \beta_2, \beta_1 \geq \alpha_1 \wedge \beta_2,
 \beta_2 \nsim \alpha_1 \vee \beta_1, \beta_2 \geq \alpha_1 \wedge \beta_1$, the following are the generators of the 
 first syzygy
 
 \begin{enumerate}
 \item $S_1=-x_{\beta_2}g_{(\alpha_1,\beta_1)} + x_{\beta_1}g_{(\alpha_1,\beta_2)}-
 x_{\alpha_1 \wedge \beta_1}g_{(\beta_2, \alpha_1 \vee \beta_1)}$.
 
 \item $S_2=x_{\alpha_1 \vee \beta_2}g_{(\alpha_1,\beta_1)}- x_{\alpha_1 \vee \beta_1}g_{(\alpha_1,\beta_2)} +
 x_{\alpha_1}g_{(\beta_2,\alpha_1 \vee \beta_1)}$.
\end{enumerate} 
\end{lemma}
\begin{proof}
 \begin{enumerate}
 \item   
$ S(f_{(\alpha_1,\beta_1)},f_{(\alpha_1,\beta_2)}) =x_{\beta_2}f_{(\alpha_1,\beta_1)}- x_{\beta_1}f_{(\alpha_1,\beta_2)}
 = x_{\beta_1}x_{\alpha_1 \vee \beta_2}x_{\alpha_1 \wedge \beta_2}-x_{\beta_2}
 x_{\alpha_1 \vee \beta_1} x_{\alpha_1 \wedge \beta_1}
 =- x_{\alpha_1 \wedge \beta_1}(x_{\beta_2}x_{\alpha_1 \vee \beta_1}- x_{\beta_1}x_{\alpha_1 \vee \beta_2})
 =-x_{\alpha_1 \wedge \beta_1}f_{(\beta_2, \alpha_1 \vee \beta_1)}
$, as $\beta_2 \nsim \alpha_1 \vee \beta_1$. \\
 Hence $-x_{\beta_2}g_{(\alpha_1,\beta_1)} + x_{\beta_1}g_{(\alpha_1,\beta_2)}-
 x_{\alpha_1 \wedge \beta_1}g_{(\beta_2, \alpha_1 \vee \beta_1)}$ is a generator of the first syzygy.

 \item 
 $ S(f_{(\alpha_1,\beta_2)},f_{(\beta_2, \alpha_1 \vee \beta_1)})= x_{\alpha_1 \vee \beta_1}f_{(\alpha_1,\beta_2)}
 - x_{\alpha_1}f_{(\beta_2, \alpha_1 \vee \beta_1)}
 =x_{\alpha_1}x_{\beta_1}x_{\alpha_1 \vee \beta_2}- x_{\alpha_1 \vee \beta_1}x_{\alpha_1 \vee \beta_2}
 x_{\alpha_1 \wedge \beta_2}
 = x_{\alpha_1 \vee \beta_2}(x_{\alpha_1}x_{\beta_1}- x_{\alpha_1 \vee \beta_1}x_{\alpha_1 \wedge \beta_1})
 =x_{\alpha_1 \vee \beta_2}f_{(\alpha_1,\beta_1)}$.\\ 
Hence, $x_{\alpha_1 \vee \beta_2}g_{(\alpha_1,\beta_1)}- x_{\alpha_1 \vee \beta_1}g_{(\alpha_1,\beta_2)} +
 x_{\alpha_1}g_{(\beta_2,\alpha_1 \vee \beta_1)}$ is a generator of the first syzygy. 
 
\end{enumerate}
\end{proof}

\begin{lemma}\label{L_1} $(\textbf{L-type})$ Let $(\alpha_1,\beta_1),(\alpha_1,\beta_2) \in \mathcal{D}_2$, for every sublattice 
 isomorphic to Figure 2.(a), with conditions $\alpha_1 \wedge \beta_1 \neq \alpha_1 \wedge \beta_2,
 \alpha_1 \vee \beta_1 \neq \alpha_1 \vee \beta_2, \beta_1 \leq \beta_2,\beta_1 \leq \alpha_1 \vee \beta_2,\beta_1 
 \nsim \alpha_1 \wedge \beta_2, \beta_2 \nsim \alpha_1 \vee \beta_1,\beta_2 \geq \alpha_1 \wedge \beta_1$ the 
 following is a generator of the first syzygy

\begin{center}
 $L=-x_{\beta_2}g_{(\alpha_1,\beta_1)} +x_{\beta_1}g_{(\alpha_1,\beta_2)} 
 + x_{\alpha_1 \vee \beta_2}g_{(\beta_1,\alpha_1 \wedge \beta_2)}
 - x_{\alpha_1 \wedge \beta_1}g_{(\beta_2, \alpha_1 \vee \beta_1)}$.
\end{center} 
\end{lemma}

\begin{proof}
$S(f_{(\alpha_1,\beta_1)},f_{(\alpha_1,\beta_2)}) =x_{\beta_2}f_{(\alpha_1,\beta_1)}- x_{\beta_1}f_{(\alpha_1,\beta_2)}
 =x_{\beta_1}x_{\alpha_1 \vee \beta_2}x_{\alpha_1 \wedge \beta_2}-x_{\beta_2}
 x_{\alpha_1 \vee \beta_1}x_{\alpha_1 \wedge \beta_1}
= -x_{\alpha_1 \wedge \beta_1}\cdot (x_{\beta_2}x_{\alpha_1 \vee \beta_1}-x_{\alpha_1 \vee \beta_2}x_{\beta_2 \wedge
(\alpha_1 \vee \beta_1)})+ x_{\alpha_1 \vee \beta_2}( x_{\beta_1}x_{\alpha_1 \wedge \beta_2}- x_{\alpha_1 \wedge \beta_1}
x_{\beta_1 \vee (\alpha_1 \wedge \beta_2)})
=-x_{\alpha_1 \wedge \beta_1} f_{(\beta_2,\alpha_1 \vee \beta_1)} + x_{\alpha_1 \vee \beta_2}\cdot
f_{(\beta_1,\alpha_1 \wedge \beta_2)} 
$,
 where the third equality follows as, 
$\beta_1 \vee (\alpha_1 \wedge \beta_2)=\beta_2 \wedge (\alpha_1 \vee \beta_1)$.\\ Hence 
$-x_{\beta_2}g_{(\alpha_1,\beta_1)} +x_{\beta_1}g_{(\alpha_1,\beta_2)} + x_{\alpha_1 \vee \beta_2}
g_{(\beta_1,\alpha_1 \wedge \beta_2)}
 - x_{\alpha_1 \wedge \beta_1}g_{(\beta_2, \alpha_1 \vee \beta_1)}$ is a generator of the first syzygy.
\end{proof} 
\begin{lemma}\label{b} $(\textbf{Box-type})$ Let $(\alpha_1,\beta_1),(\alpha_1,\beta_2) \in \mathcal{D}_2$, for every sublattice 
 isomorphic to Figure 3, with conditions $\alpha_1 \wedge \beta_1 \neq \alpha_1 \wedge \beta_2,
 \alpha_1 \vee \beta_1 \neq \alpha_1 \vee \beta_2, \beta_1 \nsim \beta_2, \beta_1 \leq \alpha_1 \vee \beta_2, 
 \beta_1 \geq \alpha_1 \wedge \beta_2, \beta_2 \nsim \alpha_1 \vee \beta_1, \beta_2 \nsim \alpha_1 \wedge \beta_1$ the 
 following are the generators of the first syzygy 

\begin{enumerate}
 \item $B_1=-x_{\beta_2}g_{(\alpha_1,\beta_1)} + x_{\beta_1}g_{(\alpha_1,\beta_2)}-x_{\alpha_1 \vee \beta_1}
 g_{(\beta_2,\alpha_1 \wedge \beta_1)}
 -x_{\alpha_1 \wedge \beta_2}g_{(\alpha_1 \vee \beta_1, \beta_1 \vee \beta_2)}$.
 
 \item $B_2=-x_{\beta_1}g_{(\alpha_1,\beta_2)}+ x_{\alpha_1}g_{(\beta_1,\beta_2)}+ x_{\beta_1 \vee \beta_2}g_{(\alpha_1,\beta_1 
 \wedge \beta_2)}+ x_{\alpha_1 \wedge \beta_2}g_{(\alpha_1 \vee \beta_1),(\beta_1 \vee \beta_2)}$.
 \end{enumerate}
\end{lemma}

\begin{proof}
 \begin{enumerate}
 \item  $S(f_{(\alpha_1,\beta_1)},f_{(\alpha_1,\beta_2)})
  =x_{\beta_2}f_{(\alpha_1,\beta_1)}- 
 x_{\beta_1}f_{(\alpha_1,\beta_2)}
 = x_{\beta_1}x_{\alpha_1 \vee \beta_2}x_{\alpha_1 \wedge \beta_2}-x_{\beta_2}
 x_{\alpha_1 \vee \beta_1} x_{\alpha_1 \wedge \beta_1}
 =-x_{\alpha_1 \vee \beta_1}(x_{\beta_2} x_{\alpha_1 \wedge \beta_1}- x_{\alpha_1 \wedge \beta_2} x_{\beta_1 \vee \beta_2})
 - x_{\alpha_1 \wedge \beta_2}(x_{\alpha_1 \vee \beta_1}x_{\beta_1 \vee \beta_2}- x_{\alpha_1 \vee \beta_2}x_{\beta_1})\\
 =-x_{\alpha_1 \vee \beta_1}f_{(\beta_2,\alpha_1 \wedge \beta_1)}- x_{\alpha_1 \wedge \beta_2}
 f_{(\alpha_1 \vee \beta_1,\beta_1 \vee \beta_2)}$.
 Hence the lemma follows.
 
\item $S(f_{(\alpha_1,\beta_2)},f_{(\beta_1,\beta_2)})
 =x_{\beta_1}f_{(\alpha_1,\beta_2)}- x_{\alpha_1}f_{(\beta_1,\beta_2)}
 =x_{\alpha_1}x_{\beta_1 \vee \beta_2}x_{\beta_1 \wedge \beta_2}- x_{\beta_1}x_{\alpha_1 \vee \beta_2}
 x_{\alpha_1 \wedge \beta_2}
 =x_{\beta_1 \vee \beta_2}\cdot(x_{\alpha_1}x_{\beta_1 \wedge \beta_2}- x_{\alpha_1 \wedge \beta_2}x_{\alpha_1 \vee \beta_1})+
 x_{\alpha_1 \wedge \beta_2}(x_{\alpha_1 \vee \beta_1} x_{\beta_1 \vee \beta_2}- x_{\beta_1}x_{\alpha_1 \vee \beta_2})
 =x_{\beta_1 \vee \beta_2}f_{(\alpha_1, \beta_1 \wedge \beta_2)} + x_{\alpha_1 \wedge \beta_2}\cdot f_{(\alpha_1 \vee \beta_1,
 \beta_1 \vee \beta_2)}$.
Hence the lemma follows.
\end{enumerate}
\end{proof}
The following theorem gives us the linear generators for the first syzygy of $R[\mathcal{L}]$.
\begin{theorem}\label{syz}Syz$^{1}_{1} R[\mathcal{L}]$ is generated by the types $S_1, S_2, L, B_1, B_2$.
\end{theorem}
\begin{proof}
 Let $f_{(\alpha_1,\beta_1)}=x_{\alpha_1}x_{\beta_1}-x_{\alpha_1 \vee \beta_1}x_{\alpha_1 \wedge \beta_1},
f_{(\alpha_2,\beta_2)}=x_{\alpha_2}x_{\beta_2}-x_{\alpha_2 \vee \beta_2}x_{\alpha_2 \wedge \beta_2}$ 
be the relation coming from the diamonds $(\alpha_1,\beta_1),(\alpha_2,\beta_2)\in \mathcal{D}_2$ . 
Therefore with respect to monomial ordering $>$, 
in$(f_{(\alpha_1,\beta_1)})=x_{\alpha_1}x_{\beta_1}$, in$(f_{(\alpha_2,\beta_2)})=x_{\alpha_2}x_{\beta_2}$. For the syzygy arising out of this S-polynomial the initial monomials must be non coprime. The only possibilities described in the Figures 2(c), 2(a) and Figure 3 and we showed in Lemmas \ref{S}, \ref{L_1}, \ref{b} that these are indeed linear syzygies.
\end{proof}
\begin{definition}\label{expressible} A pair of comparable diamonds $(\alpha_1,\beta_1), (\alpha_2,\beta_2) \in \mathcal{D}_2$ (satisfying the conditions (1) and (2)) is called \textit{expressible} if it also satisfies the condition (3) and \textit{non-expressible} otherwise.
\begin{enumerate}
  \item $\alpha_1, \beta_1, \alpha_2, \beta_2$ are all distinct.
  \item $\alpha_1 \vee \beta_1 < \alpha_2 \wedge \beta_2$.
  \item there exist $a,b \in \mathcal{L}$ such that $a\wedge b=\beta_1, a\vee b= \beta_2$ and $\alpha_1, \alpha_2, a$ each non comparable to $b$ and $\alpha_2 \wedge b$ non comparable to $\alpha_1 \vee \beta_1$.
 \end{enumerate}
 \end{definition}
 \begin{definition}\label{Diamond}
 $(\textbf{Diamond-type})$ Let Syz$_{1}^2 R[\mathcal{L}]$ be the $R[\mathcal{L}]$-module generated by  
 \begin{center}
 $D=(x_{\alpha_2}x_{\beta_2}-x_{\alpha_2 \vee \beta_2}x_{\alpha_2 \wedge \beta_2})g_{(\alpha_1,\beta_1)}-
(x_{\alpha_1}x_{\beta_1}- x_{\alpha_1 \vee \beta_1}x_{\alpha_1 \wedge \beta_1})g_{(\alpha_2,\beta_2)}$
\end{center}
for every two non-expressible diamonds $(\alpha_1,\beta_1),(\alpha_2,\beta_2)\in \mathcal{D}_2$ be called the diamond type syzygy of Syz$^{1}R[\mathcal{L}]$.
 \end{definition}

\begin{remark}
By definition of Syz$^{1}_{1}R[\mathcal{L}]$ and Syz$^{2}_{1}R[\mathcal{L}]$ we have  Syz$^{1}R[\mathcal{L}]$=Syz$^{1}_{1} R[\mathcal{L}] \oplus $Syz$^{2}_{1} R[\mathcal{L}]$. 
\end{remark}
In the following lemma, it is made clear that the $S$-polynomial arising out of two expressible diamonds is linear combination of $L$-type syzygies. Note that in the light of Lemma \ref{lem2}, for a pair of expressible diamonds $(2,3), (11,12)$ their is a diagram as in the lemma below \ref{lem3}
\begin{lemma}{\label{lem3}}In Figure 4, if $f_{(2,3)}, f_{(11,12)}$ be the relations coming from the
diamonds $(2,3),(11,12)$ respectively. Then  the $S$-polynomial arising out of these two diamonds can be expressed as a 
combination of $L$-type syzygies. 
\end{lemma}
\begin{figure}[H]
\begin{tikzpicture}[scale=.3]
 \draw (0,0)--(-2,2)--(0,4)--(-2,6)--(0,8)--(-2,10)--(0,12)--(2,10)--(4,8)--(6,6)--(4,4)--(2,2)--(0,0);
 \draw (2,2)--(0,4);
 \draw (4,4)--(2,6)--(0,8);
 \draw (0,4)--(2,6)--(4,8);
 \draw (0,8)--(2,10);
 \begin{scriptsize}
  \fill [color=black] (0,0) circle (1.5pt);
\draw[color=black] (0,-.5) node {$1$};
\fill [color=black] (-2,2) circle (1.5pt);
\draw[color=black] (-2.5,2) node {$2$};
\fill [color=black] (0,4) circle (1.5pt);
\draw[color=black] (-0.5,4) node {$4$};
\fill [color=black] (-2,6) circle (1.5pt);
\draw[color=black] (-2.5,6) node {$6$};
\fill [color=black] (0,8) circle (1.5pt);
\draw[color=black] (-.5,8) node {$9$};
\fill [color=black] (-2,10) circle (1.5pt);
\draw[color=black] (-2.5,10) node {$11$};
\fill [color=black] (0,12) circle (1.5pt);
\draw[color=black] (0,12.4) node {$13$};
\fill [color=black] (2,10) circle (1.5pt);
\draw[color=black] (2.5,10) node {$12$};
\fill [color=black] (4,8) circle (1.5pt);
\draw[color=black] (4.5,8) node {$10$};
\fill [color=black] (6,6) circle (1.5pt);
\draw[color=black] (6.5,6) node {$8$};
\fill [color=black] (4,4) circle (1.5pt);
\draw[color=black] (4.5,4) node {$5$};
\fill [color=black] (2,2) circle (1.5pt);
\draw[color=black] (2.5,2) node {$3$};
\fill [color=black] (2,6) circle (1.5pt);
\draw[color=black] (2.5,6) node {$7$};
 \end{scriptsize}
 \end{tikzpicture}
 \label{FIGURE}
 \caption{}
 \end{figure}
\begin{proof}
 Let $f_{(2,3)}=x_2x_3-x_1x_4$ and $f_{(11,12)}=x_{11}x_{12}-x_9x_{13}$ be the relations coming
 from the diamonds
 $(2,3)$ and $(11,12)$  respectively.  Then with respect to monomial ordering $>$ we have
 in$(f_{(2,3)})=x_2x_3$ and in$(f_{(11,12)})=x_{11}x_{12}$. 
 So Lcm$(f_{(2,3)},f_{(11,12)})=x_2x_3x_{11}x_{12}$.
 $S(f_{(2,3)},f_{(11,12)})=x_{11}x_{12} f_{(2,3)}-x_2x_3 f_{(11,12)}=x_2x_3x_9x_{13}-x_1x_4x_{11}x_{12}=x_9x_{13}f_1-
 x_1x_4f_6$. Therefore $(x_{11}x_{12}-x_9x_{13})g_{(2,3)}-(x_2x_3-x_1x_4)g_{(11,12)}$ 
 gives a element of the first syzygy.
 Now this expression can  be written as
 $(x_{11}x_{12}-x_9x_{13})g_{(2,3)}-(x_2x_3-x_1x_4)g_{(11,12)}= x_{11}(x_{12} g_{(2,3)}-x_6 g_{(2,8)}+ x_2 g_{(6,8)}-
 x_1 g_{(6,10)})- 
 x_{13}(x_9 g_{(2,3)}-x_6 g_{(2,5)}+ x_2 g_{(5,6)} -x_1 g_{(6,7)})+x_6(x_{13}g_{(2,5)}-x_{11}g_{(2,8)}+x_2 g_{(8,11)}-
 x_1 g_{(10,11)})+
 x_2(x_{13}g_{(5,6)}- x_{11}g_{(6,8)} +x_6 g_{(8,11)}-
 x_3 g_{(11,12)})
 -x_1(x_{13}g_{(6,7)}-x_{11}g_{(6,10)} +x_6 g_{(10,11)}- x_4 g_{(11,12)})$ which shows that it is a combination of four 
 $L$-type elements.
 Hence the lemma follows for this case.
 \end{proof}
 \begin{lemma}{\label{lem2}}
Let $(\alpha_1,\beta_1),(\alpha_2,\beta_2)\in \mathcal{D}_2$ be two expressible diamonds then there is a sublattice of $\mathcal{L}$ isomorphic to Figure 5.
\end{lemma}
\begin{figure}[H]
\begin{tikzpicture}[scale=.3]
 \draw (0,0)--(-2,2);
 \draw [dotted](0,4)--(-2,6)--(0,8);
 \draw (0,8)--(-2,10)--(0,12)--(2,10);
 \draw [dotted] (2,10)--(4,8)--(6,6)--(4,4)--(2,2);
 \draw (2,2)--(0,0);
 \draw (2,2)--(0,4);
 \draw [dotted](4,4)--(2,6)--(0,8);
 \draw [dotted](0,4)--(2,6)--(4,8);
 \draw (0,8)--(2,10);
 \begin{scriptsize}
  \fill [color=black] (0,0) circle (1.5pt);
\draw[color=black] (0,-.5) node {$\alpha_1 \wedge \beta_1$};
\fill [color=black] (-2,2) circle (1.5pt);
\draw[color=black] (-2.5,2) node {$\alpha_1$};
\fill [color=black] (0,4) circle (1.5pt);
\draw[color=black] (-1.5,4) node {$\alpha_1 \vee \beta_1$};
\fill [color=black] (-2,6) circle (1.5pt);
\draw[color=black] (-2.5,6) node {$a$};
\fill [color=black] (0,8) circle (1.5pt);
\draw[color=black] (-1.5,8) node {$\alpha_2 \wedge \beta_2$};
\fill [color=black] (-2,10) circle (1.5pt);
\draw[color=black] (-2.5,10) node {$\alpha_2$};
\fill [color=black] (0,12) circle (1.5pt);
\draw[color=black] (0,12.4) node {$\alpha_2 \vee \beta_2$};
\fill [color=black] (2,10) circle (1.5pt);
\draw[color=black] (2.7,10) node {$\beta_2$};
\fill [color=black] (4,8) circle (1.5pt);
\draw[color=black] (4.5,8) node {$z$};
\fill [color=black] (6,6) circle (1.5pt);
\draw[color=black] (6.5,6) node {$b$};
\fill [color=black] (4,4) circle (1.5pt);
\draw[color=black] (4.5,4) node {$x$};
\fill [color=black] (2,2) circle (1.5pt);
\draw[color=black] (2.7,2) node {$\beta_1$};
\fill [color=black] (2,6) circle (1.5pt);
\draw[color=black] (2.5,6) node {$y$};
 \end{scriptsize}
 \end{tikzpicture}
 \caption{}
 \label{Fig1}
 \end{figure}
\begin{proof}
In Figure 5, we have $a\wedge b=\beta_1$ and $a \vee b=\beta_2$. Now we prove that 
there exist $x,y,z \in \mathcal{L}$ such that the following holds.
\begin{enumerate}
\item There exist $x$ such that $x=\alpha_2 \wedge b=\alpha_2 \wedge \beta_2 \wedge b$ and $\beta_1 < x < b$. To prove this,
let $\alpha_2 \wedge b=x_1, (\alpha_2 \wedge \beta_2)\wedge b=x$. Now since $x_1 \vee \beta_2=
(\alpha_2 \wedge b)\vee \beta_2=\beta_2, x \vee \beta_2=((\alpha_2 \wedge \beta_2)\wedge b)\vee \beta_2=\beta_2,
x_1 \wedge \beta_2=(\alpha_2 \wedge b)\wedge \beta_2=\alpha_2 \wedge b, x \wedge \beta_2=((\alpha_2 \wedge \beta_2)\wedge b)=
\alpha_2 \wedge b$ then it implies $x=x_1$ and since we have $x=(\alpha_2 \wedge \beta_2)\wedge b$.
Let $(\alpha_2 \wedge \beta_2)\wedge b=b$. Then $\alpha_2\wedge \beta_2 \geq b$. Since $\alpha_1 > a, \beta_1 > b$, then 
 $\alpha_1\wedge \beta_1 > a\vee b= \beta_2$, which is a contradiction. Hence it follows.
 
\item $x \nsim \alpha_1 \vee \beta_1, (\alpha_1 \vee \beta_1)\wedge x=\beta_1$, as we have $x \nsim \alpha_1 \vee \beta_1, (\alpha_1 \vee \beta_1)\wedge
x=(\alpha_1 \wedge x)\vee(\beta_1\wedge x)=(\alpha_1\wedge \beta_1)\vee \beta_1=
\beta_1$.

\item $a \vee x= \alpha_2 \wedge \beta_2$ and $a \wedge x=\beta_1$, since,
$a \vee x= a \vee [(\alpha_2 \wedge \beta_2)\wedge b]=(a \vee b)\wedge [a \vee (\alpha_2 \wedge \beta_2)]
=\beta_2 \wedge (\alpha_2 \wedge \beta_2)=\alpha_2 \wedge \beta_2$ and
$a \wedge x= a \wedge [(\alpha_2 \wedge \beta_2)\wedge b]=\beta_1 \wedge (\alpha_2 \wedge \beta_2)=\beta_1$.

\item $y \wedge b=x$, where $y=(\alpha_1 \vee \beta_1)\vee x$, as
$[(\alpha_1 \vee \beta_1)\vee x]\wedge b=[(\alpha_1 \vee \beta_1)\wedge b]\vee (x \wedge b)=\beta_1 \vee x=x$.

\item $a \vee y=\alpha_2 \wedge \beta_2$ and $a \wedge y=\alpha_1 \wedge \beta_1$, since
$a \vee  [(\alpha_1 \vee \beta_1)\vee x]=[a \vee(\alpha_1 \vee \beta_1)\vee x=a \vee x= \alpha_2 \wedge \beta_2$.
and $a \wedge [(\alpha_1 \vee \beta_1)\vee x]=(a \wedge x)\vee [a \wedge(\alpha_1 \vee \beta_1)]
=\beta_1 \vee (\alpha_1 \vee \beta_1 )=\alpha_1 \vee \beta_1$.

\item $(\alpha_2 \wedge \beta_2)\wedge z=y$ and $(\alpha_2 \wedge \beta_2)\vee z=\beta_2$, where 
$z=(\alpha_1 \vee \beta_1)\vee b$, as
 $(\alpha_2 \wedge \beta_2)\wedge [(\alpha_1 \vee \beta_1)\vee b]
=[(\alpha_2 \wedge \beta_2)\wedge (\alpha_1 \vee \beta_1)]\vee [(\alpha_2 \wedge \beta_2)\wedge b]
=(\alpha_1 \vee \beta_1)\vee x=y$
and $(\alpha_2 \wedge \beta_2)\vee [(\alpha_1 \vee \beta_1)\vee b]
=(\alpha_2 \wedge \beta_2)\vee (\alpha_1 \vee \beta_1)=(\alpha_2 \wedge \beta_2)\vee b=\beta_2$.
 
\end{enumerate}
\end{proof}

\begin{lemma}\label{expressible_main_lemma} For every pair of diamonds $(\alpha_1,\beta_1), (\alpha_2,\beta_2)\in \mathcal{D}_2$ the syzygy
 \begin{center}
 $Z=(x_{\alpha_2}x_{\beta_2}-x_{\alpha_2 \vee \beta_2}x_{\alpha_2 \wedge \beta_2})g_{(\alpha_1,\beta_1)}-
(x_{\alpha_1}x_{\beta_1}- x_{\alpha_1 \vee \beta_1}x_{\alpha_1 \wedge \beta_1})g_{(\alpha_2,\beta_2)}$
\end{center}
is congruent to a syzygy $S$ in Syz$_{1}^{2}R[\mathcal{L}]$ modulo Syz$_{1}^{1}R[\mathcal{L}]$. Further if a pair of comparable diamonds are non-expressible then the syzygy element $S$ is not zero. 
\end{lemma}
\begin{proof}
 If the syzygy $Z$ is expressible then by Lemma \ref{lem3} it belongs to Syz$_{1}^{1}R[\mathcal{L}]$. If $Z$ is not expressible then there will not exist a sublattice as in Figure 5. (see Lemma \ref{lem2}). So, without loss of generality one of the following cases will arise.
 \vskip 1mm 
 \noindent \textsc{Case 1:} $\alpha_1 \leq \alpha_2, \alpha_1 \leq \beta_2, \beta_1 \leq \alpha_2, \beta_1 \nsim \beta_2$\\
 In this case, we can have the lattices isomorphic to the following sublattices Figure 6. (a) and (b).
  \begin{figure}[H] 
 \begin{tikzpicture}[scale=.3]
 \draw (0,0)--(-2,2)--(0,4)--(-2,6)--(0,8)--(2,10)--(4,8)--(2,6)--(4,4)--(2,2)--(0,0);
 \draw (2,2)--(0,4)--(2,6)--(0,8);
 \draw (10,0)--(8,2)--(10,4)--(12,6)--(14,4)--(12,2)--(10,0);
 \draw (12,2)--(10,4);
 \draw (20,0)--(18,2)--(16,4)--(16,4)--(18,6)--(20,8)--(22,6)--(24,4)--(22,2)--(20,0);
 \draw (18,2)--(20,4)--(22,6);
 \draw (22,2)--(20,4)--(18,6);
 \draw (32,10)--(34,8)--(36,6)--(34,4)--(32,2)--(30,0)--(28,2)--(26,4)--(28,6)--(30,8)--(32,10);
 \draw (28,2)--(30,4)--(32,6)--(34,8);
 \draw (32,2)--(30,4)--(28,6);
 \draw (34,4)--(32,6)--(30,8);
 
 \draw (-1,-1.5) node[anchor=north west] {$(a)$};
 \draw (9,-1.5) node[anchor=north west] {$(b)$};
 \draw (19,-1.5) node[anchor=north west] {$(c)$};
 \draw (29,-1.5) node[anchor=north west] {$(d)$};
 \begin{scriptsize}
  \fill [color=black] (0,0) circle (1.5pt);
  \draw[color=black] (0,-.5) node {$1$};
  \fill [color=black] (-2,2) circle (1.5pt);
  \draw[color=black] (-3.5,2) node {$\alpha_1=2$};
  \fill [color=black] (0,4) circle (1.5pt);
  \draw[color=black] (-.5,4) node {$4$};
\fill [color=black] (-2,6) circle (1.5pt);
\draw[color=black] (-3.5,6) node {$\alpha_2=6$};
\fill [color=black] (0,8) circle (1.5pt);
\draw[color=black] (-.5,8) node {$8$};
\fill [color=black] (2,10) circle (1.5pt);
\draw[color=black] (2,10.5) node {$10$};
\fill [color=black] (4,8) circle (1.5pt);
\draw[color=black] (5.5,8) node {$9=\beta_2$};
\fill [color=black] (2,6) circle (1.5pt);
\draw[color=black] (2.5,6) node {$7$};
\fill [color=black] (4,4) circle (1.5pt);
\draw[color=black] (5.5,4) node {$5=\beta_1$};
\fill [color=black] (2,2) circle (1.5pt);
\draw[color=black] (2.5,2) node {$3$};
\fill [color=black] (10,0) circle (1.5pt);
\fill [color=black] (8,2) circle (1.5pt);
\draw[color=black] (7.5,2) node {$\beta_1$};
\fill [color=black] (10,4) circle (1.5pt);
\draw[color=black] (9.4,4) node {$\alpha_2$};
\fill [color=black] (12,6) circle (1.5pt);
\fill [color=black] (14,4) circle (1.5pt);
\draw[color=black] (14.5,4) node {$\beta_2$};
\fill [color=black] (12,2) circle (1.5pt);
\draw[color=black] (12.6,2) node {$\alpha_1$};
\fill [color=black] (20,0) circle (1.5pt);
\fill [color=black] (18,2) circle (1.5pt);
\fill [color=black] (16,4) circle (1.5pt);
\draw[color=black] (15.5,4) node {$\alpha_1$};
\fill [color=black] (18,6) circle (1.5pt);
\fill [color=black] (20,8) circle (1.5pt);
\fill [color=black] (22,6) circle (1.5pt);
\fill [color=black] (24,4) circle (1.5pt);
\draw[color=black] (24.5,4) node {$\beta_2$};
\fill [color=black] (22,2) circle (1.5pt);
\draw[color=black] (22.6,2) node {$\beta_1$};
\fill [color=black] (20,4) circle (1.5pt);
\draw[color=black] (19.5,4) node {$\alpha_2$};
\fill [color=black] (30,0) circle (1.5pt);
\fill [color=black] (28,2) circle (1.5pt);
\fill [color=black] (26,4) circle (1.5pt);
\draw[color=black] (25.5,4) node {$\alpha_1$};
\fill [color=black] (28,6) circle (1.5pt);
\draw[color=black] (27.4,6) node {$\alpha_2$};
\fill [color=black] (30,8) circle (1.5pt);
\fill [color=black] (32,10) circle (1.5pt);
\fill [color=black] (30,4) circle (1.5pt);
\fill [color=black] (32,6) circle (1.5pt);
\draw[color=black] (32.6,6) node {$\beta_1$};
\fill [color=black] (32,2) circle (1.5pt);
\fill [color=black] (34,4) circle (1.5pt);
\fill [color=black] (34,8) circle (1.5pt);
\fill [color=black] (36,6) circle (1.5pt);
\draw[color=black] (36.6,6) node {$\beta_2$};
 \end{scriptsize}
\label{fig}
 \end{tikzpicture}
 \caption{}
   \end{figure}
 \noindent In Figure 6. (a), numbering the lattice we have the syzygy from the diamonds $(\alpha_1,\beta_1), (\alpha_2,\beta_2)$ that is $(2,5), (6,9)$ is 
 
 $Z=(x_6x_9-x_4x_{10})g_{(2,5)}-(x_2x_5-x_1x_7)g_{(6,9)}=x_2[-x_9g_{(5,6)}+x_5g_{(6,9)}+x_{10}g_{(4,5)}-x_3g_{(8,9)}]+x_9[-x_6g_{(2,5)}+x_2g_{(5,6)}+x_8g_{(2,3)}-x_1g_{(6,7)}]+x_{10}[x_4g_{(2,5)}-x_2g_{(4,5)}-x_7g_{(2,3)}]+x_1[x_4g_{(8,9)}-x_7g_{(6,9)}+x_9g_{(6,7)}]+(x_2x_3-x_1x_4)g_{(8,9)}-(x_8x_9-x_7x_{10})g_{(2,3)}$,
 which is $S=(x_2x_3-x_1x_4)g_{(8,9)}-(x_8x_9-x_7x_{10})g_{(2,3)}$ modulo Syz$_{1}^{1}R[\mathcal{L}]$ and it belongs to Syz$_{1}^{2}R[\mathcal{L}]$.
 In Figure 6. (b), the lattice is Grid lattice $G(1,2)$, then the syzygy $Z$ is in Syz$_{1}^{1}R[\mathcal{L}]$.
 \vskip 1mm \noindent
 \textsc{Case 2:} $\alpha_1 \nsim \alpha_2, \alpha_1 \nsim \beta_2, \beta_1 \leq \alpha_2, \beta_1 \nsim \beta_2$\\
 In this case, we can have a sublattice isomorphic to the lattice Figure 5.(c).
 Since the lattice Figure 6.(c), is Grid lattice $G(2,2)$ then the syzygy $Z$ belongs to Syz$_{1}^{1}R[\mathcal{L}]$.
 \vskip 1mm \noindent
 \textsc{Case 3:} $\alpha_1 \leq \alpha_2, \alpha_1 \nsim \beta_2, \beta_1 \nsim \alpha_2, \beta_1 \nsim \beta_2$\\
 In this case, we can have a sublattice isomorphic to the following lattice Figure 6.(d).
 Since the lattice is Grid lattice $G(2,3)$ then it belongs to Syz$_{1}^{1}R[\mathcal{L}]$.
Therefore, by above we have for a pair of comparable diamonds which are non-expressible is not zero. Hence the claim.
\end{proof}

\begin{corollary}\label{linear_indep_diamond}
Let $(D_{i_1},D_{j_1}), (D_{i_2},D_{j_2}),\ldots, (D_{i_n},D_{j_n})$ be $n$-pairs of distinct non-expressible diamonds and let $f_1, f_2,\ldots,f_n$ be the corresponding syzygies. Then $f_1, f_2,\ldots, f_n$ are $k$-linearly independent.
\end{corollary}
\begin{proof}
The syzygy element $f_i$ has the expression $f_i = (x_{\alpha_2}x_{\beta_2}-x_{\alpha_2 \vee \beta_2}x_{\alpha_2 \wedge \beta_2})g_{(\alpha_1,\beta_1)}-
(x_{\alpha_1}x_{\beta_1}- x_{\alpha_1 \vee \beta_1}x_{\alpha_1 \wedge \beta_1})g_{(\alpha_2,\beta_2)}$, where  $(\alpha_1,\beta_1), (\alpha_2,\beta_2)\in \mathcal{D}_2$ are the two diamonds. Since the two diamonds are distinct, then $g_{(\alpha_i, \beta_i)}$ are distinct basis elements for $i \neq j$. So, clearly these are $k$-linearly independent. 
\end{proof}
\begin{remark}
Using the corollary \ref{linear_indep_diamond} we can find the rank of the module Syz$_1^2(\mathcal{L})$, we call it by the number $n_{\mathcal{D}}R[\mathcal{L}]$.
\end{remark}
\begin{lemma}
\label{lemma_reverse}
  Let $(\alpha_1,\beta_1),(\alpha_2,\beta_2)\in \mathcal{D}_2$, the syzygy $Z=(x_{\alpha_2}x_{\beta_2}-x_{\alpha_2 \vee \beta_2}x_{\alpha_2 \wedge \beta_2})g_{(\alpha_1,\beta_1)}-
(x_{\alpha_1}x_{\beta_1}- x_{\alpha_1 \vee \beta_1}x_{\alpha_1 \wedge \beta_1})g_{(\alpha_2,\beta_2)}$ arising out of these two diamonds is
  expressible by strip, $L$, box types then there is a sublattice as in Figure 7.
 \end{lemma}
 \begin{figure}[H]
\begin{tikzpicture}[scale=.35]
\draw (0,0)--(-1,1)--(0,2)--(-1,3)--(-2,4)--(-3,5)--(-2,6)--(-1,7)--(0,8)--(-1,9)--(0,10)--(1,9)--(2,8)--(3,7)--(4,6)--(5,5)--(4,4)--(3,3)--(2,2)--(1,1)--(0,0);
\draw (1,1)--(0,2);
\draw (2,2)--(1,3)--(0,4)--(-1,5)--(-2,6);
\draw (4,4)--(3,5)--(2,6)--(1,7)--(0,8);
\draw (0,2)--(1,3)--(2,4)--(3,5)--(4,6);
\draw (-2,4)--(-1,5)--(0,6)--(1,7)--(2,8);
 \draw (0,8)--(1,9);
 \draw [dashed] (-1,3)--(0,4)--(1,5)--(2,6)--(3,7);
 \draw [dashed] (3,3)--(2,4)--(1,5)--(0,6)--(-1,7);
 \begin{scriptsize}
\fill [color=black] (0,0) circle (1.5pt);
\fill [color=black] (-1,1) circle (1.5pt);
\fill [color=black] (0,2) circle (1.5pt);
\fill [color=black] (-1,3) circle (1.5pt);
\fill [color=black] (-2,4) circle (1.5pt);
\fill [color=black] (-3,5) circle (1.5pt);
\fill [color=black] (-2,6) circle (1.5pt);
\fill [color=black] (-1,7) circle (1.5pt);
\fill [color=black] (0,8) circle (1.5pt);
\fill [color=black] (-1,9) circle (1.5pt);
\fill [color=black] (0,10) circle (1.5pt);
\fill [color=black] (1,9) circle (1.5pt);
\fill [color=black] (2,8) circle (1.5pt);
\fill [color=black] (3,7) circle (1.5pt);
\fill [color=black] (4,6) circle (1.5pt);
\fill [color=black] (5,5) circle (1.5pt);
\fill [color=black] (4,4) circle (1.5pt);
\fill [color=black] (3,3) circle (1.5pt);
\fill [color=black] (2,2) circle (1.5pt);
\fill [color=black] (1,1) circle (1.5pt);
\fill [color=black] (1,3) circle (1.5pt);
\fill [color=black] (0,4) circle (1.5pt);
\fill [color=black] (-1,5) circle (1.5pt);
\fill [color=black] (2,4) circle (1.5pt);
\fill [color=black] (1,5) circle (1.5pt);
\fill [color=black] (0,6) circle (1.5pt);
\fill [color=black] (3,5) circle (1.5pt);
\fill [color=black] (2,6) circle (1.5pt);
\fill [color=black] (1,7) circle (1.5pt);
\draw[color=black] (-1.7,1) node {$\alpha_1$};
\draw[color=black] (1.7,1) node {$\beta_1$};
\draw[color=black] (-1.7,9) node {$\alpha_2$};
\draw[color=black] (1.7,9) node {$\beta_2$};
\draw[color=black] (-2,5) node {$D_1$};
\draw[color=black] (4,5) node {$D_2$};
 \end{scriptsize}
 \end{tikzpicture}
 \label{fig1}
 \caption{}
 \end{figure}
 
 \begin{proof}
We will show that if there is no sublattice isomorphic to the Figure 6, then the syzygy element
 $(x_{\alpha_2}x_{\beta_2}-x_{\alpha_2 \vee \beta_2}x_{\alpha_2 \wedge \beta_2})g_{(\alpha_1,\beta_1)}-
 (x_{\alpha_1}x_{\beta_1}-x_{\alpha_1 \vee \beta_1}x_{\alpha_1 \wedge \beta_1})g_{(\alpha_2,\beta_2)}$ cannot be written as a combination of the linear syzygies. We argue in the following way if there is no sublattice isomorphic to the one in the figure, it must be a sublattice of the same. We will show for maximal sublaticces of the syzygy element cannot be written as a combination of linear syzygies thus arriving at a contradiction.
 
 A maximal sublattice must have either one of the two corner diamonds are missing. So let either the diamond $D_1$ or $D_2$ (in Figure 7) be missing.  Let us say $D_1$ is missing.
 Then $x_{\alpha_2}g_{(\alpha_1,\beta_1)}$ or $x_{\beta_2}g_{(\alpha_1,\beta_1)}$ is a term of a syzygy generator. But from the 
 figure we see that the coefficient of $g_{(\alpha_1,\beta_1)}$ could only be the points which are thick in the figure
 and we see that $x_{\alpha_2}$ or $x_{\beta_2}$ is not there. Thus this element cannot be expressed as a combination of
 other type of syzygy generators. Similarly, we can show that if the diamond $D_2$ is not there.
 \end{proof}

\begin{theorem}
\label{main}
 Syz$^{2}_{1} R[\mathcal{L}]$ is trivial if and only if every pair of comparable diamonds in $\mathcal{D}_2$ are expressible.
 \end{theorem}
  \begin{proof}
  If every pair of diamonds in the lattice $\mathcal{L}$ is expressible then by definition the module Syz$_1^2R[\mathcal{L}]$ is trivial. Now, for the other direction let us assume the module Syz$_1^2R[\mathcal{L}]$ is trivial and further let us assume that there is a syzygy  $Z$ which is not expressible. As a result of the Lemma \ref{expressible_main_lemma} we know that the syzygy element associated to this syzygy is in Syz$_1^2R[\mathcal{L}]$ modulo Syz$_1^1R[\mathcal{L}]$. Now, since Syz$_1^2R[\mathcal{L}]$ is trivial it must be congruent to zero modulo the linear syzygy module. Thus the syzygy $Z$ is expressible. 
  \end{proof}
\begin{remark}\label{Grid-lattice}
The Theorem \ref{main}, tells us when the first syzygy is linear, i.e. when all the comparable diamonds are expressible. In particular, for the grid lattice $G(m,n)$ the first syzygies are linear. 
\end{remark}
In the following section we will refine our result for linearity of the first syzygy even stronger by showing that it only depends on the number of pairs of join-meet irreducibles. Let $JM$ be the set of join-meet irreducibles of the lattice $\mathcal{L}$ and let $k(\mathcal{L})=\# \{(\theta_i,\theta_j): \theta_i \nsim \theta_j, \theta_i,\theta_j \in JM\}$.

\begin{theorem}
 \label{last} With above notations we have the following
\begin{enumerate}
 \item When $k(\mathcal{L})=1$, the first syzygy of $R[\mathcal{L}]$ is linear.
 \item When $k(\mathcal{L})=2$,
 \begin{enumerate}
  \item If the lattice is Figure 8, then first syzygy is non linear.
  \item  Else if the lattice is Figure 9, then first syzygy is linear if and only if $ht(\theta_2 \wedge \theta_3)
 -ht(\theta_1 \wedge \theta_2)=1$ or $ht(\theta_2 \vee \theta_3)-ht(\theta_1 \vee \theta_2)=1$, where $\theta_1,\theta_2,
 \theta_3 \in JM$.
 \end{enumerate}
 \item If $k(\mathcal{L}) \geq 3$ then the first syzygy is non linear.
\end{enumerate} 
\end{theorem}
\begin{proof}
 \begin{enumerate}
 \item Since $k(\mathcal{L})=1$ then it is a grid lattice. Therefore, by Remark \ref{Grid-lattice} it is linear.
 \begin{figure}[H]
\begin{tikzpicture}[scale=.3]
 \draw (0,0)--(-2,2)--(0,4)--(-2,6)--(0,8)--(2,6)--(0,4)--(2,2)--(0,0);
 \draw[dashed] (.8,.8)--(-1.2,2.8);
 \draw[dashed] (1.2,1.2)--(-.8,3.2);
 \draw[dashed] (-.8,.8)--(1.2,2.8);
 \draw[dashed] (-1.2,1.2)--(.8,3.2);
 \draw[dashed] (-.8,4.8)--(1.2,6.8);
 \draw[dashed] (-1.2,5.12)--(.8,7.12);
 \draw[dashed] (.8,4.8)--(-1.2,6.8);
 \draw[dashed] (1.2,5.12)--(-.8,7.12);

 \begin{scriptsize}
  \fill [color=black] (0,0) circle (1.5pt);
\fill [color=black] (-2,2) circle (1.5pt);
\draw[color=black] (-2.5,2) node {$\theta_1$};
\fill [color=black] (0,4) circle (1.5pt);
\fill [color=black] (-2,6) circle (1.5pt);
\draw[color=black] (-2.5,6) node {$\theta_3$};
\fill [color=black] (0,8) circle (1.5pt);
\fill [color=black] (2,6) circle (1.5pt);
\draw[color=black] (2.5,6) node {$\theta_4$};
\fill [color=black] (2,2) circle (1.5pt);
\draw[color=black] (2.5,2) node {$\theta_2$};
 \end{scriptsize}
 \end{tikzpicture}
 \label{diamond}
 \caption{}
 \end{figure}
 
 \item 
 \begin{enumerate}  
  \item If the lattice as in Figure 8, then the first syzygy have diamond type generators such as 
  $(x_{\theta_3}x_{\theta_4}-x_{\theta_3 \vee \theta_4}x_{\theta_3 \wedge \theta_4})g_{(\theta_1,\theta_2)}
  -(x_{\theta_1}x_{\theta_2}-x_{\theta_1 \vee \theta_2}x_{\theta_1 \wedge \theta_2})g_{(\theta_3,\theta_4)}$ and it cannot be 
  expressed by linear generators. Therefore, it belongs to Syz$_{1}^{2}R[\mathcal{L}]$ by Lemma \ref{expressible_main_lemma}. Hence it is non linear.
    \begin{figure}[H]
\begin{tikzpicture}[scale=.3]
 \draw (0,0)--(-2,2)--(0,4)--(-2,6)--(0,8)--(2,6)--(0,4)--(2,2)--(0,0);
 \draw (2,2)--(4,4)--(2,6);
 \draw[dashed] (.8,.8)--(-1.2,2.8);
 \draw[dashed] (1.2,1.2)--(-.8,3.2);
 \draw[dashed] (-.8,.8)--(3.4,4.8);
 \draw[dashed] (-1.2,1.2)--(3.2,5.2);
 \draw[dashed] (-.8,4.8)--(1.2,6.8);
 \draw[dashed] (-1.2,5.12)--(.8,7.12);
 \draw[dashed] (2.8,2.8)--(-1.2,6.8);
 \draw[dashed] (3.12,3.12)--(-.8,7.12);

 \begin{scriptsize}
  \fill [color=black] (0,0) circle (1.5pt);
  \draw[color=black] (0,-.4) node {$\theta_1 \wedge \theta_2$};
\fill [color=black] (-2,2) circle (1.5pt);
\draw[color=black] (-2.5,2) node {$\theta_1$};
\fill [color=black] (0,4) circle (1.5pt);
\fill [color=black] (-2,6) circle (1.5pt);
\draw[color=black] (-2.5,6) node {$\theta_3$};
\fill [color=black] (0,8) circle (1.5pt);
\draw[color=black] (0,8.3) node {$\theta_2 \vee \theta_3$};
\fill [color=black] (2,6) circle (1.5pt);
\draw[color=black] (3.5,6) node {$\theta_1 \vee \theta_2$};
\fill [color=black] (2,2) circle (1.5pt);
\draw[color=black] (3.5,2) node {$\theta_2 \wedge \theta_3$};
\fill [color=black] (4,4) circle (1.5pt);
\draw[color=black] (4.5,4) node {$\theta_2$};
 \end{scriptsize}
 \end{tikzpicture}
 \label{Ltype}
 \caption{}
 \end{figure}

   \item $\Rightarrow$: If the lattice as in Figure 9, and since first syzygy is linear then there is a sublattice as in (see Lemma \ref{lem2}) and therefore $ht(\theta_2 \wedge \theta_3) -ht(\theta_1 \wedge \theta_2)=1$
   or $ht(\theta_2 \vee \theta_3)-ht(\theta_1 \vee \theta_2)=1$, where $\theta_1,\theta_2, \theta_3 \in JM$.\\   
   $\Leftarrow$: Let the condition holds then all diamond type elements of the first syzygy can be expressed by other types.
   Hence the first syzygy is linear.
 
   \end{enumerate}
  \begin{figure}[H]
\begin{tikzpicture}[scale=.3]
 \draw (0,0)--(-2,2)--(0,4)--(-2,6)--(0,8)--(2,6)--(0,4)--(2,2)--(0,0);
 \draw (2,2)--(4,4)--(2,6);
 \draw (0,0)--(2,-2)--(4,0)--(2,2);
 \draw[dashed] (.8,.8)--(-1.2,2.8);
 \draw[dashed] (1.2,1.2)--(-.8,3.2);
 \draw[dashed] (-.8,.8)--(3.4,4.8);
 \draw[dashed] (-1.2,1.2)--(3.2,5.2);
 \draw[dashed] (-.8,4.8)--(1.2,6.8);
 \draw[dashed] (-1.2,5.12)--(.8,7.12);
 \draw[dashed] (2.8,2.8)--(-1.2,6.8);
 \draw[dashed] (3.12,3.12)--(-.8,7.12);
 \draw[dashed] (.8,.8)--(2.8,-1.2);
  \draw[dashed] (1.2,1.2)--(3.2,-1.);
  \draw[dashed] (.8,-.8)--(2.8,1.2);
   \draw[dashed] (1.2,-1.2)--(3.12,.8);
 \begin{scriptsize}
  \fill [color=black] (0,0) circle (1.5pt);
  \fill [color=black] (-2,2) circle (1.5pt);
\draw[color=black] (-2.5,2) node {$\theta_2$};
\fill [color=black] (0,4) circle (1.5pt);
\fill [color=black] (-2,6) circle (1.5pt);
\draw[color=black] (-2.5,6) node {$\theta_4$};
\fill [color=black] (0,8) circle (1.5pt);
\fill [color=black] (2,6) circle (1.5pt);
\fill [color=black] (2,2) circle (1.5pt);
\fill [color=black] (4,4) circle (1.5pt);
\draw[color=black] (4.5,4) node {$\theta_3$};
\fill [color=black] (4,0) circle (1.5pt);
\draw[color=black] (4.5,0) node {$\theta_1$};
 \end{scriptsize}
 \end{tikzpicture}
 \label{L2}
 \caption{}
 \end{figure}
   
\item If $k(\mathcal{L}) \geq 3$ then we have a figure as in Figure 10, then the first syzygy element coming from the 
 lower diamond and
 upper diamond cannot be expressed by other type of syzygies. Hence diamond type elements will occur. Thus for this case the 
 first syzygy is not linear.
\end{enumerate}
\end{proof}
\section{First Betti number for planar distributive lattice}\label{Betti}

 Let $\mathcal{L}$ be a planar distributive lattice and $JM$ be the set of all join-meet irreducible elements of $\mathcal{L}$.
Let $D_{(\theta_i,\theta_j)}=[\theta_i \wedge \theta_j,\theta_i \vee \theta_j]$, the interval for $\theta_i \nsim \theta_j, 
\theta_i, \theta_j \in JM$.
One can write the lattice $\mathcal{L}$ as a union of these intervals i.e. union of grid lattices. Let $n(\theta_1,\theta_2)$ be the number of strip 
type generator in $D_{(\theta_1,\theta_2)}$. Then we have the following 
\begin{lemma}
For $\theta_1,\theta_2 \in JM$,
\begin{align*}
 n(\theta_1,\theta_2)=& S(ht(\theta_1)-ht(\theta_1 \wedge \theta_2), ht(\theta_2)-ht(\theta_1 \wedge \theta_2))
\end{align*}
where $S(.,.)$ is the number of strip type generator from the grid lattice.
\end{lemma}
\begin{proof} Since the grid lattice formed by $D_{(\theta_1,\theta_2)}$ is the lattice 
$G(ht(\theta_1)-ht(\theta_1 \wedge \theta_2),ht(\theta_2)-ht(\theta_1 \wedge \theta_2))$. Hence it follows.
\end{proof}
\begin{lemma}
\label{one}
 For a planar distributive lattice $\mathcal{L}$, the number of strip type generators of first syzygy of $R[\mathcal{L}]$, we denote it by 
 $n_{S}(\mathcal{L})$ is the following
 \begin{align*}
 n_{S}(\mathcal{L}) =& \sum_{\substack{\theta_i \nsim \theta_j \\ \theta_i,\theta_j \in JM}} n(\theta_i,\theta_j)
  - \sum_{\theta_i,\theta_j,\theta_k} n(\theta_i \vee(\theta_j \wedge \theta_k),\theta_j).
\end{align*}
 \end{lemma}
 
 \begin{figure}[H] 
 \begin{tikzpicture}[scale=.3]
 \draw (0,0)--(-2,2)--(-4,4)--(-2,6)--(0,8)--(2,6)--(0,4)--(2,2)--(0,0);
 \draw (-2,2)--(0,4);
 \draw (0,4)--(-2,6);
 
 \begin{scriptsize}
 \fill [color=black] (0,0) circle (1.5pt);
\draw[color=black] (0,-.5) node {$\theta_i\wedge \theta_j$};
\fill [color=black] (-2,2) circle (1.5pt);
\draw[color=black] (-3.4,2) node {$\theta_j \wedge \theta_k$};
\fill [color=black] (-4,4) circle (1.5pt);
\draw[color=black] (-4.5,4) node {$\theta_j$};
\fill [color=black] (0,4) circle (1.5pt);
\draw[color=black] (2.7,4) node {$\theta_i \vee (\theta_j \wedge \theta_k)$};
\fill [color=black] (2,2) circle (1.5pt);
\draw[color=black] (2.5,2) node {$\theta_i$};
\fill [color=black] (-2,6) circle (1.5pt);
\draw[color=black] (-3.4,6) node {$\theta_i \vee \theta_j$};
\fill [color=black] (0,8) circle (1.5pt);
\draw[color=black] (0,8.4) node {$\theta_j \vee \theta_k$};
\fill [color=black] (2,6) circle (1.5pt);
\draw[color=black] (2.5,6) node {$\theta_k$};

 \end{scriptsize}
\end{tikzpicture}
\label{L}
\caption{}
\end{figure}
 \begin{proof}
From the Figure 11. it follows.
\end{proof}
Let $m(\theta_1,\theta_2)$ be the number of $L$-type generator in $D_{(\theta_1,\theta_2)}$. Then we have the following
\begin{lemma}
 For $\theta_1,\theta_2 \in JM$,
 \begin{align*}
  m(\theta_1,\theta_2)=& L(ht(\theta_1)-ht(\theta_1 \wedge \theta_2),ht(\theta_2)-ht(\theta_1 \wedge \theta_2)).
 \end{align*}
where $L(.,.)$ denotes the number of $L$-shape generator from the grid lattice.
\end{lemma}
\begin{proof}
Since, the number of $L$-type generators in $D_{(\theta_1,\theta_2)}$ is same as to find the number of $L$-type generators in grid lattice, then it follows.
\end{proof}
\begin{lemma}
\label{two}
 For a planar distributive lattice $\mathcal{L}$, the number of $L$-type generator of $R[\mathcal{L}]$, we denote it by 
 $n_{L}(\mathcal{L})$ is the following
 \begin{align*} 
  n_{L}(\mathcal{L})=& \sum_{\substack{\theta_i \nsim \theta_j \\ \theta_1,\theta_j \in JM}}m(\theta_i,\theta_j) -\sum L (r,s)\\&
  + \sum_{\theta_i,\theta_j,\theta_k}(ht(\theta_j \vee \theta_k)-ht(\theta_i \vee \theta_j)
  (ht(\theta_j \wedge \theta_k)-ht(\theta_i \wedge \theta_j)) 
  \binom{r+1}{2} \binom{s+1}{2}.
   \end{align*}
   where $r=ht(\theta_i \vee \theta_j)-ht(\theta_j), s=ht(\theta_j)-ht(\theta_j \wedge \theta_k)$.
\end{lemma}

\begin{figure}[H]
\begin{tikzpicture}[scale=.35]
 \draw (0,0)--(-2,2)--(0,4)--(-2,6)--(0,8)--(2,6)--(0,4)--(2,2)--(0,0);
 \draw (2,2)--(4,4)--(2,6);
 \draw[dashed] (0,1.8)--(2.2,4)--(-0.4,6.6)--(-0.8,6.2)--(1.5,4)--(-.4,2.2)--(0,1.8);

 \begin{scriptsize}
  \fill [color=black] (0,0) circle (1.5pt);
  \draw[color=black] (0,-.4) node {$\theta_i \wedge \theta_j$};
\fill [color=black] (-2,2) circle (1.5pt);
\draw[color=black] (-2.5,2) node {$\theta_i$};
\fill [color=black] (0,4) circle (1.5pt);
\fill [color=black] (-2,6) circle (1.5pt);
\draw[color=black] (-2.5,6) node {$\theta_k$};
\fill [color=black] (0,8) circle (1.5pt);
\draw[color=black] (0,8.3) node {$\theta_j \vee \theta_k$};
\fill [color=black] (2,6) circle (1.5pt);
\draw[color=black] (3.5,6) node {$\theta_i \vee \theta_j$};
\fill [color=black] (2,2) circle (1.5pt);
\draw[color=black] (3.5,2) node {$\theta_j \wedge \theta_k$};
\fill [color=black] (4,4) circle (1.5pt);
\draw[color=black] (4.5,4) node {$\theta_j$};
 \end{scriptsize}
 \end{tikzpicture}
 \label{dot}
 \caption{}
 \end{figure}

\begin{proof}
The number of $L$-type for each $D_{(\theta_i,\theta_j)}$ is $m(\theta_i,\theta_j)$. For each $\theta_i,
\theta_j,\theta_k$ there is double counting $m(r,s)$ where $r=ht(\theta_i \vee \theta_j)-ht(\theta_j),
s=ht(\theta_j)-ht(\theta_j \wedge \theta_k)$. Also for each $\theta_i,\theta_j,\theta_k$ there is a $L$-type (see the dotted 
line in the the Figure 11) and the number is $(ht(\theta_j \vee \theta_k)-ht(\theta_i \vee \theta_j)
  (ht(\theta_j \wedge \theta_k)-ht(\theta_i \wedge \theta_j)) \binom{r+1}{2} \binom{s+1}{2}$. Since for $(\theta_i,\theta_j),
 (\theta_j,\theta_k),(\theta_k,\theta_l)$ there is no $L$-type. Hence the lemma follows.
  \end{proof}

\noindent Let $B(\theta_i,\theta_j)$ be the number of box type in $D_{(\theta_1,\theta_2)}$ then we have the following:

\begin{lemma}
\label{three}
 For a planar distributive lattice $\mathcal{L}$, the number of box type generator of $R[\mathcal{L}]$, $n_{B}(\mathcal{L})$ is
 given by
 \begin{align*}
  n_{B}(\mathcal{L})=& \sum_{\substack{\theta_i \nsim \theta_j \\ \theta_1,\theta_j \in JM}}B(\theta_i,\theta_j)
  -\sum_{\theta_i,\theta_j,\theta_k}B(\theta_i \vee (\theta_j \wedge \theta_k),\theta_j).
 \end{align*}
\end{lemma}
\begin{proof}
For each $\theta_i,\theta_j$ the number of box types is $B(\theta_i,\theta_j)$ and for each $\theta_i,
\theta_j,\theta_k$ there is double counting $B(\theta_i \vee (\theta_j \wedge \theta_k),\theta_j)$ in numbers in their
intersections.
Hence the lemma follows.
\end{proof}
\begin{theorem}\label{betti} The first Betti number $\beta_1(\mathcal{L})$ of $R[\mathcal{L}]$ is 
\begin{align*}
 \beta_1(\mathcal{L})=&n_{S}(\mathcal{L})+n_{L}(\mathcal{L})+n_{B}(\mathcal{L})+n_{D}(\mathcal{L}).
\end{align*}
\end{theorem}
\begin{proof}
It follows from the Lemmas \ref{one}, \ref{two}, \ref{three} and Corollary \ref{linear_indep_diamond}.
\end{proof}

In the following example the results of this paper is demonstrated. 
\begin{example} 
Let $\mathcal{L}$ be a lattice as shown in the following Figure 13.

\begin{minipage}{\linewidth}
\begin{minipage}{0.3\linewidth}
\begin{figure}[H]
\begin{tikzpicture}[scale=.3]
 \draw (0,0)--(-2,2)--(0,4)--(-2,6)--(0,8)--(2,10)--(0,12)--(2,14)--(4,12)--(6,10)--(4,8)--(6,6)--(4,4)--(2,2)--(0,0);
 \draw (2,2)--(0,4)--(2,6)--(4,8);
 \draw (4,4)--(2,6)--(0,8);
 \draw (4,8)--(2,10)--(4,12);
 \begin{scriptsize}
 \fill [color=black] (0,0) circle (1.5pt);
\draw[color=black] (0,-.5) node {$1$};
\fill [color=black] (-2,2) circle (1.5pt);
\draw[color=black] (-3.4,2) node {$\theta_1=2$};
\fill [color=black] (2,2) circle (1.5pt);
\draw[color=black] (2.5,2) node {$3$};
\fill [color=black] (0,4) circle (1.5pt);
\draw[color=black] (-.7,4) node {$4$};
\fill [color=black] (4,4) circle (1.5pt);
\draw[color=black] (4.7,4) node {$5$};
\fill [color=black] (-2,6) circle (1.5pt);
\draw[color=black] (-3.4,6) node {$\theta_2=6$};
\fill [color=black] (2,6) circle (1.5pt);
\draw[color=black] (2.7,6) node {$7$};
\fill [color=black] (6,6) circle (1.5pt);
\draw[color=black] (7.5,6) node {$8=\theta_3$};
\fill [color=black] (0,8) circle (1.5pt);
\draw[color=black] (-.5,8) node {$9$};
\fill [color=black] (4,8) circle (1.5pt);
\draw[color=black] (4.7,8) node {$10$};
\fill [color=black] (2,10) circle (1.5pt);
\draw[color=black] (1,10) node {$11$};
\fill [color=black] (6,10) circle (1.5pt);
\draw[color=black] (7.5,10) node {$12=\theta_4$};
\fill [color=black] (0,12) circle (1.5pt);
\draw[color=black] (-1.5,12) node {$\theta_5=13$};
\fill [color=black] (4,12) circle (1.5pt);
\draw[color=black] (4.7,12) node {$14$};
\fill [color=black] (2,14) circle (1.5pt);
\draw[color=black] (2.7,14) node {$15$};
 \end{scriptsize}
\end{tikzpicture}
\label{Ex}
\caption{}
\end{figure}
\end{minipage}
\begin{minipage}{0.6\linewidth}Here the set of join-meet irreducible elements is $JM=\{ \theta_1,\theta_2,\theta_3,\theta_4,\theta_5\}$. Since cardinality of $k(\mathcal{L}) \geq 3$, so the first syzygy of $R[\mathcal{L}]$ is non linear follows from the Theorem \ref{last}.
In $\mathcal{L}$, we can see that there is a sublattice isomorphic to the Figure 14.(a), so from this we get two strip type syzygies $S_1=-x_5 g_{(2,3)}+x_3 g_{(2,5)}-x_1 g_{(4,5)}$ and $S_2= x_7 g_{(2,3)}-x_4 g_{(2,5)}+x_2 g_{(4,5)}$. Since there is a sublattice isomorphic to the Figure 14.(b), there is a $L$-type syzygy $L=-x_{10} g_{(6,10)}+ x_5 g_{(6,5)}+x_{11}g_{(5,4)}-x_3 g_{(10,9)}$. Also, since in $\mathcal{L}$, there is a sublattice isomorphic to Figure 14.(c), one can get two box type syzygy $B_1=-x_8 g_{(6,7)}+x_7 g_{(6,8)}-x_9g_{(8,4)}-x_3g_{(9,10)}$ and $B_2=-x_7 g_{(6,8)}+x_6 g_{(7,8)}+x_{10}g_{(6,5)}+x_3 g_{(9,10)}$.
\end{minipage}
\end{minipage}

\begin{figure}[H]
\begin{tikzpicture}[scale=.25]
 \draw (0,0)--(-2,2)--(0,4)--(2,6)--(4,4)--(2,2)--(0,0);
 \draw (2,2)--(0,4);
 \draw (12,0)--(10,2)--(8,4)--(10,6)--(12,8)--(14,6)--(12,4)--(14,2)--(12,0);
 \draw (10,2)--(12,4)--(10,6);
 \draw (22,0)--(20,2)--(18,4)--(20,6)--(22,8)--(24,6)--(26,4)--(24,2)--(22,0);
 \draw (20,2)--(22,4)--(24,6);
 \draw (24,2)--(22,4)--(20,6);
 \draw (-.5,-1) node[anchor=north west] {$(a)$};
 \draw (11.5,-1) node[anchor=north west] {$(b)$};
 \draw (21.5,-1) node[anchor=north west] {$(c)$};
 \begin{scriptsize}
 \fill [color=black] (0,0) circle (1.5pt);
 \draw[color=black] (0,-.5) node {$1$};
\fill [color=black] (-2,2) circle (1.5pt);
\draw[color=black] (-2.7,2) node {$2$};
\fill [color=black] (0,4) circle (1.5pt);
\draw[color=black] (-.7,4.5) node {$4$};
\fill [color=black] (2,6) circle (1.5pt);
\draw[color=black] (2,6.5) node {$7$};
\fill [color=black] (4,4) circle (1.5pt);
\draw[color=black] (4.5,4) node {$5$};
\fill [color=black] (2,2) circle (1.5pt);
\draw[color=black] (2.7,2) node {$3$};
\fill [color=black] (12,0) circle (1.5pt);
\draw[color=black] (12,-.5) node {$3$};
\fill [color=black] (10,2) circle (1.5pt);
\draw[color=black] (9,2) node {$4$};
\fill [color=black] (8,4) circle (1.5pt);
\draw[color=black] (7.3,4) node {$6$};
\fill [color=black] (10,6) circle (1.5pt);
\draw[color=black] (9.5,6.5) node {$9$};
\fill [color=black] (12,8) circle (1.5pt);
\draw[color=black] (12,8.5) node {$11$};
\fill [color=black] (14,6) circle (1.5pt);
\draw[color=black] (14.7,6) node {$10$};
\fill [color=black] (12,4) circle (1.5pt);
\draw[color=black] (12.7,4) node {$7$};
\fill [color=black] (14,2) circle (1.5pt);
\draw[color=black] (14.5,2) node {$5$};
\fill [color=black] (22,0) circle (1.5pt);
\draw[color=black] (22,-.5) node {$3$};
\fill [color=black] (20,2) circle (1.5pt);
\draw[color=black] (19.3,2) node {$4$};
\fill [color=black] (18,4) circle (1.5pt);
\draw[color=black] (17.3,4) node {$6$};
\fill [color=black] (20,6) circle (1.5pt);
\draw[color=black] (19.5,6.5) node {$9$};
\fill [color=black] (22,8) circle (1.5pt);
\draw[color=black] (22,8.5) node {$11$};
\fill [color=black] (22,4) circle (1.5pt);
\draw[color=black] (22.7,4) node {$7$};
\fill [color=black] (24,6) circle (1.5pt);
\draw[color=black] (24.7,6.3) node {$10$};
\fill [color=black] (26,4) circle (1.5pt);
\draw[color=black] (26.7,4) node {$8$};
\fill [color=black] (24,2) circle (1.5pt);
\draw[color=black] (24.7,2) node {$5$};
 \end{scriptsize}
\end{tikzpicture}
\label{figure}
\caption{}
\end{figure}
 By Definition \ref{expressible}, there is a diamond type coming from the two diamonds $(2,3)$ and $(13,14)$ and the two diamonds $(2,3)$ and $(9,10)$ are expressible by linear syzygies.
\end{example}

\section{$\beta_1$ of Hibi ideals with regularity at most 3}

In the following, as an application of our results especially \ref{betti}, we will give closed form formulae for the first Betti numbers of Hibi ideals with a given regularity $r\leq 3$. At this juncture let us recall the results from \cite{herzog}, in which the authors gave a complete characterization of the regularity of a distributive lattice $\mathcal{L}$, namely $\mbox{reg} I(\mathcal{L})= |P| - \mbox{rank} (P)$, where $P$ is the poset of nontrivial join-irreducible elements of $\mathcal{L}$. In the next subsection we will calculate the number of generators of strip, L and Box type syzygies for the ``grid lattice", which will be used repeatedly while applying Lemmas \ref{one}, \ref{two}, \ref{three} and finally in form of the Theorem \ref{betti}. We also use the Corollary \ref{linear_indep_diamond} to calculate the number of linearly independent quadratic generators. Putting these together, we will give closed form formulae for all the planar lattices for regularity less than or equal to three.
\subsection{The first Betti number of Grid lattice}\label{grid_betti}

{\bf Strip type generators:} Let $n_S(G(1,n))$ be the number of strip type generators in a grid lattice $G(1,n)$, we have 
\begin{center}
	$n_S(G(1,n))=2n_S(G(1,n-1))- n_S(G(1,n-2) + 2(n-1)=2 \binom{n+1}{3}.$
\end{center}

\noindent Let $n_S(G(m,n))$ be the number of strip type generators for the grid lattice $G(m,n)$ we have,
\begin{align*}
	n_S(G(m,n))&=[m+(m-1)+ \ldots +1] n_S(G(1,n))+[n+(n-1)+ \ldots + 1] n_S(G(1,m))\\
	&=\binom{m+1}{2}n_S(G(1,n))+\binom{n+1}{2}n_S(G(1,m)).
\end{align*}
{\bf $L$-type generators:} Let $n_L(G(2,n))$ be the number of $L$-type generators for grid lattice $G(2,n)$ we have,
\begin{align*}
	n_L(G(2,n))&=2 \sum_{i=1}^{n}(n+1-i)(i-1)\\
	&=\frac{n(n^2-1)}{3}.
\end{align*}
Therefore, the $L$-type generators say $n_L(G(m,n))$, for grid lattice $G(m,n)$ is,
\begin{align*}
	n_L(G(m,n))& = \sum_{i=0}^{m}(m+1-i)(i-1)n_L(G(2,n))\\
	&=\frac{m(m^2-1)}{6}n_L(G(2,n))\\
	&= \frac{n_L(G(2,m))n_L(G(2,n))}{2}.
\end{align*}
{\bf Box type generators:} Since in a grid lattice clearly the number of box type generators is same as the number of $L$-type generators, hence, if $n_B(G(m,n))$ is the number of box type generators of $G(m,n)$ we have,
\begin{center}
	$n_B(G(m,n))=\frac{n_B(G(2,m))n_B(G(2,n))}{2}$
\end{center}
Therefore, $\beta_1$ of $R[\mathcal{L}]$ for $G(m,n)$ is 
\begin{center}
	$\beta_1(\mathcal{L})= n_S(G(m,n))+ n_L(G(m,n))+ n_B(G(m,n))$.
	\end{center}
\subsection{The first Betti number of planar lattices of regularity 2} In the following, lefthand side figure represent join irreducible poset and righthand side figure represent the corresponding lattice. Figure $15$ and $16$ are the regularity $2$ lattices up to isomorphism.


	\begin{minipage}{\linewidth}
		\begin{minipage}{0.3\linewidth}
			\begin{figure}[H]
				\begin{tikzpicture}[scale=.15]
			\draw (-2,0)--(-2,2);
			\draw [dotted](-2,4)--(-2,2);
			\draw (-2,4)--(-2,6);
			\draw (6,0)--(4,2)--(6,4)--(8,6)--(10,4)--(8,2)--(6,0);
			\draw (10,8)--(12,10)--(14,8)--(12,6)--(10,8);
			\draw (6,4)--(8,2);
			\draw [dotted] (8,6)--(10,8);
			\draw [dotted](10,4)--(12,6);
			\draw [decorate,decoration={brace,amplitude=5pt,mirror,raise=.5pt}, yshift=0pt](6.5,-.5)--(14.5,7.5) node [black,midway,xshift=0.5 cm]{\footnotesize $n$};
			\begin{scriptsize}
			\fill [color=black] (-5,0) circle (4 pt);
			\fill [color=black] (-2,0) circle (4 pt);
			\fill [color=black] (-2,2) circle (4 pt);
			\fill [color=black] (-2,4) circle (4 pt);
			\fill [color=black] (-2,6) circle (4 pt);
			\fill [color=black] (6,0) circle (4 pt);
			\fill [color=black] (4,2) circle (4 pt);
			\fill [color=black] (8,2) circle (4 pt);
			\fill [color=black] (6,4) circle (4 pt);
			\fill [color=black] (10,4) circle (4 pt);
			\fill [color=black] (8,6) circle (4 pt);
			\fill [color=black] (12,6) circle (4 pt);
			\fill [color=black] (10,8) circle (4 pt);
			\fill [color=black] (14,8) circle (4 pt);
			\fill [color=black] (12,10) circle (4 pt);
			\end{scriptsize}
				\end{tikzpicture}
				\caption{}
			\end{figure}
		\end{minipage}
		\begin{minipage}{0.6\linewidth}Since here the lattice is $G(1,n)$, therefore, the first Betti number is same as the first Betti number of $G(m,n)$, for $m=1$, (see \ref{grid_betti}). We can see that in this lattice it has no $L$-type, box type and diamond type. Hence all numbers are the following
			\begin{align*}
				n_{S}(\mathcal{L}) &=2\binom{n+1}{3}\\
				n_{L}(\mathcal{L}) &=0\\
			   n_{B}(\mathcal{L}) &=0\\
				n_{D}(\mathcal{L}) &=0
			\end{align*}
			Therefore first Betti number $\beta_1(\mathcal{L})=2 \binom{n+1}{3}$.
		\end{minipage}
	\end{minipage}

	\begin{minipage}{\linewidth}
		\begin{minipage}{0.3\linewidth}
			\begin{figure}[H]
				\begin{tikzpicture}[scale=.15]
				\draw (-5,0)--(-2,3);
				\draw (-2,0)--(-2,2);
				\draw [dotted](-2,4)--(-2,2);
				\draw (-2,4)--(-2,6);
				\draw (6,0)--(4,2)--(6,4)--(8,2)--(6,0);
				\draw (6,4)--(8,6);
				\draw [dotted] (8,6)--(10,8);
				\draw (10,8)--(12,10);
				
				\begin{scriptsize}
				\fill [color=black] (-5,0) circle (4 pt);
				\fill [color=black] (-2,0) circle (4 pt);
				\fill [color=black] (-2,2) circle (4 pt);
				\fill [color=black] (-2,4) circle (4 pt);
				\fill [color=black] (-2,6) circle (4 pt);
				\fill [color=black] (6,0) circle (4 pt);
				\fill [color=black] (4,2) circle (4 pt);
				\fill [color=black] (8,2) circle (4 pt);
				\fill [color=black] (6,4) circle (4 pt);
				\fill [color=black] (8,6) circle (4 pt);
				\fill [color=black] (10,8) circle (4 pt);
				\fill [color=black] (12,10) circle (4 pt);
				\end{scriptsize}
				\end{tikzpicture}
				\caption{}
			\end{figure}
		\end{minipage}
		\begin{minipage}{0.6\linewidth}
					Since here there is no strip, $L$, box and diamond types, therefore,
				$n_{S}(\mathcal{L})=n_{L}(\mathcal{L})=n_{B}(\mathcal{L})=n_{D}(\mathcal{L})=0$.
			As a result $\beta_1(\mathcal{L})=0$.
		\end{minipage}
	\end{minipage}

\subsection{The first Betti number of planar lattices of regularity 3}
The following are the list of regularity $3$ lattices up to isomorphic.

	\begin{minipage}{\linewidth}
		\begin{minipage}{0.3\linewidth}
			\begin{figure}[H]
				\begin{tikzpicture}[scale=.15]
				\draw (-5,0)--(-5,2);
				\draw (-2,0)--(-2,2);
				\draw (-2,4)--(-2,6);
	            \draw [decorate,decoration={brace,amplitude=5pt,mirror,raise=.5pt}, yshift=0pt](8.5,-.5)--(16.5,7.5) node [black,midway,xshift=0.5 cm]{\footnotesize $n$};		
				\draw [dotted] (-2,2)--(-2,4);
				\draw (8,0)--(6,2)--(4,4)--(6,6)--(8,8)--(10,6)--(12,4)--(10,2)--(8,0);
				\draw (6,2)--(8,4)--(10,6);
				\draw [dotted] (8,8)--(10,10);
				\draw [dotted] (12,4)--(14,6);
				\draw [dotted] (10,6)--(12,8);
				\draw (14,6)--(16,8)--(14,10)--(12,12)--(10,10)--(12,8)--(14,10);
				\draw (6,6)--(8,4)--(10,2);
				\draw (14,6)--(12,8);	
				
				\begin{scriptsize}
				
				\fill [color=black] (-5,0) circle (4 pt);
				\fill [color=black] (-5,2) circle (4 pt);
				\fill [color=black] (-2,0) circle (4 pt);
				\fill [color=black] (-2,2) circle (4 pt);
				\fill [color=black] (-2,4) circle (4 pt);
				\fill [color=black] (-2,6) circle (4 pt);
				\fill [color=black] (8,0) circle (4 pt);
				\fill [color=black] (6,2) circle (4 pt);
				\fill [color=black] (4,4) circle (4 pt);
				\fill [color=black] (6,6) circle (4 pt);
				\fill [color=black] (8,8) circle (4 pt);
				\fill [color=black] (10,10) circle (4 pt);
				\fill [color=black] (12,12) circle (4 pt);
				\fill [color=black] (14,10) circle (4 pt);
				\fill [color=black] (16,8) circle (4 pt);
				\fill [color=black] (14,6) circle (4 pt);
				\fill [color=black] (12,8) circle (4 pt);
				\fill [color=black] (10,6) circle (4 pt);
				\fill [color=black] (8,4) circle (4 pt);
				\fill [color=black] (10,2) circle (4 pt);
				\fill [color=black] (12,4) circle (4 pt);
				\end{scriptsize}
				\end{tikzpicture}
				\caption{}
				\label{F1}
			\end{figure}
		\end{minipage}
		\begin{minipage}{0.6\linewidth} The figure \ref{F1} is $G(2,n)$, therefore $\beta_1$ is same as $\beta_1$ of $G(m,n)$ for $m=2$ (see \ref{grid_betti}). Here we do not have any diamond types, therefore all the numbers are as below,
			\begin{align*}
				n_{S}(\mathcal{L}) &=n^2(n+1)\\
				n_{L}(\mathcal{L}) &=\frac{n(n^2-1)}{3}\\
				n_{B}(\mathcal{L}) &=\frac{n(n^2-1)}{3}\\
				n_{D}(\mathcal{L}) &=0
			\end{align*}
			Therefore, $\beta_1(\mathcal{L})=\frac{n(5n^2+3n-2)}{3}$.
		\end{minipage}
	\end{minipage}
	\\ \\
 \begin{figure}[H]
			\begin{tikzpicture}[scale=.15]
			\draw (-5,0)--(-5,4);
			\draw (-2,0)--(-2,2);
			\draw (-2,4)--(-2,6);
			\draw (-5,0)--(-2,3);
			\draw [dotted] (-2,2)--(-2,4);
			\draw (8,0)--(6,2)--(4,4)--(6,6)--(8,4)--(10,2)--(8,0);
			\draw [dotted] (6,6)--(8,8);
			\draw [dotted] (8,4)--(10,6);
			\draw [dotted] (10,2)--(12,4);
			\draw (8,8)--(10,10)--(12,12)--(14,10)--(12,8)--(14,6)--(12,4);
			\draw [dotted] (14,10)--(16,12);
			\draw [dotted] (12,12)--(14,14);
			\draw  (14,14)--(16,16)--(18,14)--(16,12)--(14,14);
			\draw (8,8)--(10,6)--(12,4);
			\draw (10,10)--(12,8);
			\draw (6,2)--(8,4);
			\draw (10,6)--(12,8);
			 \draw [decorate,decoration={brace,amplitude=5pt,mirror,raise=.5pt}, yshift=0pt](8.5,-.5)--(14.5,5.5) node [black,midway,xshift=0.5 cm]{\footnotesize $m$};
			  \draw [decorate,decoration={brace,amplitude=5pt},xshift=-4pt, yshift=0pt](3.5,4.5)--(15.5,16.5) node [black,midway,xshift=-.5 cm]{\footnotesize $n$};
			\draw (-.5,-1) node[anchor=north west] {$(a)$};
			\draw (33.5,-1) node[anchor=north west] {$(b)$};
		\draw (30,0)--(30,4);
			\draw (33,0)--(33,2);
			\draw (33,4)--(33,6);
			\draw (30,4)--(33,3);
			\draw [dotted] (33,2)--(33,4);
			\draw (42,0)--(40,2)--(42,4)--(44,2)--(42,0);
			\draw [dotted] (42,4)--(44,6);
			\draw [dotted] (44,2)--(46,4);
			\draw (46,4)--(44,6)--(46,8)--(48,10)--(46,12)--(48,10)--(50,8)--(48,6)--(46,4);
			\draw (48,6)--(46,8);
			\draw [dotted] (46,12)--(48,14);
			\draw [dotted] (48,10)--(50,12);
			\draw [dotted] (50,8)--(52,10);
			\draw (52,10)--(50,12)--(48,14)--(50,16)--(52,14)--(54,12)--(52,10);
			\draw (50,12)--(52,14);
			\draw (46,8)--(48,10);
			\draw (46,8)--(44,10)--(46,12);
			 \draw [decorate,decoration={brace,amplitude=5pt,mirror,raise=.5pt}, yshift=0pt](42.5,-.5)--(54.5,11.5) node [black,midway,xshift=0.5 cm]{\footnotesize $m$};
			\draw [decorate,decoration={brace,amplitude=5pt},xshift=-4pt, yshift=0pt](43.5,10.5)--(49.5,16.5) node [black,midway,xshift=-.5 cm]{\footnotesize $n$};
			\begin{scriptsize}
			\fill [color=black] (-5,0) circle (4 pt);
			\fill [color=black] (-5,4) circle (4 pt);
			\fill [color=black] (-2,0) circle (4 pt);
			\fill [color=black] (-2,2) circle (4 pt);
			\fill [color=black] (-2,4) circle (4 pt);
			\fill [color=black] (-2,6) circle (4 pt);
			\fill [color=black] (8,0) circle (4 pt);
			\fill [color=black] (6,2) circle (4 pt);
			\fill [color=black] (4,4) circle (4 pt);
			\fill [color=black] (10,2) circle (4 pt);
			\fill [color=black] (8,4) circle (4 pt);
			\fill [color=black] (6,6) circle (4 pt);
			\fill [color=black] (12,4) circle (4 pt);
			\fill [color=black] (10,6) circle (4 pt);
			\fill [color=black] (8,8) circle (4 pt);
			\fill [color=black] (14,6) circle (4 pt);
			\fill [color=black] (12,8) circle (4 pt);
			\fill [color=black] (10,10) circle (4 pt);
			\fill [color=black] (14,10) circle (4 pt);
			\fill [color=black] (12,12) circle (4 pt);
			\fill [color=black] (16,12) circle (4 pt);
			\fill [color=black] (14,14) circle (4 pt);
			\fill [color=black] (16,16) circle (4 pt);
			\fill [color=black] (18,14) circle (4 pt);
			\fill [color=black] (30,0) circle (4 pt);
			\fill [color=black] (30,4) circle (4 pt);
			\fill [color=black] (33,0) circle (4 pt);
			\fill [color=black] (33,2) circle (4 pt);
			\fill [color=black] (33,4) circle (4 pt);
			\fill [color=black] (33,6) circle (4 pt);
			\fill [color=black] (42,0) circle (4 pt);
			\fill [color=black] (40,2) circle (4 pt);
			\fill [color=black] (44,2) circle (4 pt);
			\fill [color=black] (42,4) circle (4 pt);
			\fill [color=black] (46,4) circle (4 pt);
			\fill [color=black] (44,6) circle (4 pt);
			\fill [color=black] (48,6) circle (4 pt);
			\fill [color=black] (46,8) circle (4 pt);
			\fill [color=black] (44,10) circle (4 pt);
			\fill [color=black] (50,8) circle (4 pt);
			\fill [color=black] (48,10) circle (4 pt);
			\fill [color=black] (46,12) circle (4 pt);
			\fill [color=black] (52,10) circle (4 pt);
			\fill [color=black] (50,12) circle (4 pt);
			\fill [color=black] (48,14) circle (4 pt);
			\fill [color=black] (54,12) circle (4 pt);
			\fill [color=black] (52,14) circle (4 pt);
			\fill [color=black] (50,16) circle (4 pt);
			\end{scriptsize}
			\end{tikzpicture}
			\caption{}
			\label{F2}
		\end{figure}
The Figure \ref{F2} $(a)$ (and by symmetry Figure \ref{F2}, $(b)$) is a union of two grid lattices $G(2,m)$ and $G(1,n)$ with an ``overlap" of $G(1,m)$. Therefore, by the  Lemmas \ref{one}, \ref{two}, \ref{three}, the numbers are as below, note that there is no diamond types, and $\beta_1$ is sum of these,
		\begin{align*}
		n_{S}(\mathcal{L}) &=n_{S}(G(2,m))+ n_S(G(1,n))- n_S(G(1,m))\\
		 n_{L}(\mathcal{L}) &=
		 \begin{cases}
                  n_{L}(G(2,m))+(n-m) \binom{m+1}{2} & m \geq 2	\\
			n-1 & m=1\\
			\end{cases}
			\\
		n_{B}(\mathcal{L}) &=
			\begin{cases}
			\frac{m(m^2-1)}{3} &  m \geq 2	\\
			0 &  m=1\\
			\end{cases}
			\\
	n_{D}(\mathcal{L}) &=0
		\end{align*}

In figure \ref{F4}, it is clear that, the number of strip type generators is sum of the strip types in $G(m,1)$ and $G(1,n)$ of the two ``legs" of this figure. $L$-types are particular cases of $(2)$ and $(3)$. We can see that there is no box types. The diamond types are obtained by coupling any diamond from the upper leg with these of the lower leg (avoiding the elbow) by virtue of \ref{linear_indep_diamond}. Consolidating all that, we have,
\begin{minipage}{\linewidth}
	\begin{minipage}{0.3\linewidth}
		\begin{figure}[H]
			\begin{tikzpicture}[scale=.15]
			\draw (-5,0)--(-5,2);
			\draw (-5,0)--(-2,2)--(-2,0);
			\draw (-2,2)--(-2,4);
			\draw (-2,8)--(-2,6);
			\draw (-5,-4)--(-5,-2);
			\draw [dotted] (-2,6)--(-2,4);
			\draw [dotted] (-5,0)--(-5,-2);
			\draw (12,-4)--(10,-2)--(12,0)--(14,-2)--(12,-4);
			\draw [dotted] (10,-2)--(8,0);
			\draw [dotted] (12,0)--(10,2);
			\draw (8,0)--(6,2)--(4,4)--(6,6)--(8,8)--(10,6)--(8,4)--(10,2)--(8,0);
			\draw (6,6)--(8,4)--(6,2);
			\draw [dotted](8,8)--(10,10);
			\draw [dotted](10,6)--(12,8);
			\draw (12,8)--(10,10)--(12,12)--(14,10)--(12,8);
			 \draw [decorate,decoration={brace,amplitude=5pt}, xshift=-4pt,yshift=0pt](11.5,-4.5)--(5.5,1.5) node [black,midway,xshift=-0.5 cm]{\footnotesize $m$};
			\draw [decorate,decoration={brace,amplitude=5pt},xshift=-4pt, yshift=0pt](5.5,6.5)--(11.5,12.5) node [black,midway,xshift=-.5 cm]{\footnotesize $n$};
			
			\begin{scriptsize}
			\fill [color=black] (-5,0) circle (4 pt);
			\fill [color=black] (-5,2) circle (4 pt);
			\fill [color=black] (-2,0) circle (4 pt);
			\fill [color=black] (-2,2) circle (4 pt);
			\fill [color=black] (-2,4) circle (4 pt);
			\fill [color=black] (-2,6) circle (4 pt);
			\fill [color=black] (-2,8) circle (4 pt);
			\fill [color=black] (-5,-2) circle (4 pt);
			\fill [color=black] (-5,-4) circle (4 pt);
			\fill [color=black] (12,-4) circle (4 pt);
			\fill [color=black] (14,-2) circle (4 pt);
			\fill [color=black] (10,-2) circle (4 pt);
			\fill [color=black] (12,0) circle (4 pt);
			\fill [color=black] (8,0) circle (4 pt);
			\fill [color=black] (10,2) circle (4 pt);
			\fill [color=black] (6,2) circle (4 pt);
			\fill [color=black] (8,4) circle (4 pt);
			\fill [color=black] (4,4) circle (4 pt);
			\fill [color=black] (6,6) circle (4 pt);
			\fill [color=black] (8,8) circle (4 pt);
			\fill [color=black] (10,6) circle (4 pt);
			\fill [color=black] (10,10) circle (4 pt);
			\fill [color=black] (12,8) circle (4 pt);
			\fill [color=black] (12,12) circle (4 pt);
			\fill [color=black] (14,10) circle (4 pt);				
			\end{scriptsize}
			\end{tikzpicture}
			\caption{}
			\label{F4}
		\end{figure}
	\end{minipage}
	\begin{minipage}{0.6\linewidth} 
		\begin{align*}
			n_{S}(\mathcal{L}) &=2\Big[\binom{m+2}{3}+\binom{n+2}{3}\Big] \\
			n_{L}(\mathcal{L}) &=mn\\
			n_{B}(\mathcal{L}) &=0\\						
			n_{D}(\mathcal{L}) &=
			\begin{cases}
			\binom{m}{2}\binom{n}{2} & m,n \geq 2	\\
			0 &  m=1 \, \text{or}\, n=1\\
			\end{cases}
		\end{align*}
		Therefore, $\beta_1(\mathcal{L})$ is sum of above.
	\end{minipage}
\end{minipage}

\begin{minipage}{\linewidth}
		\begin{minipage}{0.3\linewidth}
		\begin{figure}[H]
			\begin{tikzpicture}[scale=.15]
			\draw (-5,0)--(-5,6);
		\draw (-2,0)--(-2,2);
		\draw (-2,6)--(-2,8);
		\draw (-5,6)--(-2,3);
		\draw (-5,0)--(-2,4.5);
		\draw [dotted] (-2,2)--(-2,6);
		\draw (6,0)--(4,2)--(6,4)--(8,2)--(6,0);
		\draw [dotted] (6,4)--(8,6);
		\draw [dotted] (8,2)--(10,4);
		\draw (12,10)--(10,8)--(12,6)--(10,4);
		\draw (10,8)--(8,10)--(10,12);
		\draw [dotted](10,12)--(12,14);
		\draw [dotted](12,10)--(14,12);
		\draw (12,14)--(14,16)--(16,18)--(18,16)--(16,14)--(18,12)--(16,10)--(14,12)--(12,14);
		\draw (14,8)--(12,10)--(10,12);
		\draw (16,10)--(14,12)--(12,14);
		\draw (12,6)--(14,8);
		\draw [dotted](14,8)--(16,10);
		\draw (10,4)--(8,6)--(10,8);
		\draw (14,12)--(16,14)--(14,16);
		\draw [dotted](16,18)--(18,20);
		\draw (18,20)--(20,22)--(22,20)--(20,18)--(18,20);
		\draw [dotted](18,16)--(20,18);	
		 \draw [decorate,decoration={brace,amplitude=5pt,mirror,raise=.5pt}, yshift=0pt](6.5,-.5)--(12.5,5.5) node [black,midway,xshift=0.5 cm]{\footnotesize $m$};
		\draw [decorate,decoration={brace,amplitude=5pt},xshift=-4pt, yshift=0pt](7.5,10.5)--(13.5,16.5) node [black,midway,xshift=-.5 cm]{\footnotesize $k$};
		\draw [decorate,decoration={brace,amplitude=5pt},xshift=-4pt, yshift=0pt](13.5,16.5)--(19.5,22.5) node [black,midway,xshift=-.5 cm]{\footnotesize $n$};	
		
		\begin{scriptsize}
		\fill [color=black] (-5,0) circle (4 pt);
		\fill [color=black] (-5,6) circle (4 pt);
		\fill [color=black] (-2,0) circle (4 pt);
		\fill [color=black] (-2,2) circle (4 pt);
		\fill [color=black] (-2,6) circle (4 pt);
		\fill [color=black] (-2,8) circle (4 pt);
		\fill [color=black] (6,0) circle (4 pt);
		\fill [color=black] (4,2) circle (4 pt);
		\fill [color=black] (8,2) circle (4 pt);
		\fill [color=black] (6,4) circle (4 pt);
		\fill [color=black] (10,4) circle (4 pt);
		\fill [color=black] (8,6) circle (4 pt);
		\fill [color=black] (12,6) circle (4 pt);
		\fill [color=black] (10,8) circle (4 pt);
		\fill [color=black] (8,10) circle (4 pt);
		\fill [color=black] (14,8) circle (4 pt);
		\fill [color=black] (12,10) circle (4 pt);
		\fill [color=black] (10,12) circle (4 pt);
		\fill [color=black] (16,10) circle (4 pt);
		\fill [color=black] (14,12) circle (4 pt);
		\fill [color=black] (12,14) circle (4 pt);
		\fill [color=black] (18,12) circle (4 pt);
		\fill [color=black] (16,14) circle (4 pt);
		\fill [color=black] (14,16) circle (4 pt);
		\fill [color=black] (18,16) circle (4 pt);
		\fill [color=black] (16,18) circle (4 pt);
		\fill [color=black] (20,18) circle (4 pt);
		\fill [color=black] (18,20) circle (4 pt);
		\fill [color=black] (22,20) circle (4 pt);
		\fill [color=black] (20,22) circle (4 pt);
		\end{scriptsize}
		\end{tikzpicture}
			\caption{}
			\label{F5}
		\end{figure}
		\end{minipage}
		\begin{minipage}{0.6\linewidth} In the figure \ref{F5}, note that there are strip types coming from the two legs and the middle portion (a copy of $G(2,k)$), as a result we use our lemma \ref{one} to obtain the number of strip types. The number of box types are coming from the copy of $G(2,k)$, and the $L$ types are computed using lemma \ref{two}. The diamond types can be obtained by coupling all the diamonds from the upper leg with that of lower leg but we have to remove the ones that are expressible as a result of the lemma \ref{linear_indep_diamond}.
		\end{minipage}
	\end{minipage}
	\begin{align*}
			n_{S}(\mathcal{L}) &=n_S(G(1,m+k))+ n_S(G(1,n+k))+ n_S(G(2,k))-2 n_S(G(1,k)),  \, k\geq 0\\
			n_{L}(\mathcal{L})& =
			\begin{cases}
			n_{L}(G(2,k))+(n+m) \binom{k+1}{2} & k \geq 2	\\
			mn &  k=1\\
			\end{cases}
	\\
			n_{B}(\mathcal{L}) &=
			\begin{cases}
		      \frac{k(k^2-1)}{3}   & k \geq 2	\\
		     	0 & k< 2\\
			\end{cases}
		\\
			n_{D}(\mathcal{L}) &=(D(m)+mk)(D(n)+nk)-D(k)mn,\,  \text{where} \,D(m)=\frac{m(m+1)}{2}
		\end{align*}
		Therefore, $\beta_1(\mathcal{L})$ is sum of above numbers.

\textbf{Acknowledgement:} The authors would like to thank the referee for many helpful comments. The first author would like to thank BITS Pilani K.K. Birla Goa campus, Goa, India.

\bibliographystyle{abbrv}
\bibliography{refs_syzygy}
\end{document}